\documentclass[11pt]{article}
\usepackage{latexsym}
\usepackage{amsfonts,amssymb,amsmath,mathtools}
\usepackage{bbold}

\newtheorem{thm}{Theorem}[section]
\newtheorem{cor}[thm]{Corollary}
\newtheorem{lem}[thm]{Lemma}
\newtheorem{pro}[thm]{Proposition}
\newtheorem{defn}[thm]{Definition}

\newcommand{\minus}{\smallsetminus}
\newcommand{\x}{\hspace{-0.025in}\times\hspace{-0.025in}}
\newcommand{\e}{\varepsilon}

\title{The word problem of the Brin-Higman-Thompson groups}

\author{J.C.\ Birget}

\date{\scriptsize
26 vi 2020}

\begin{document}
\maketitle

\begin{abstract}
We show that the word problem of the Brin-Higman-Thompson group $n G_{k,1}$ 
is {\sf coNP}-complete for all $n \ge 2$ and all $k \ge 2$.
For this we prove that $n G_{k,1}$ is finitely generated, and that 
$n G_{k,1}$ contains a subgroup of $2 G_{2,1}$ that can represent bijective
circuits. 

We also show that for all $n \ge 1$ and
$k \ge 2$: \ If $\,K = 1 + (k-1)\,N\,$ for some $N \ge 1$, then
$n G_{K,1} \le n G_{k,1}$. In particular, $n G_{K,1} \le n G_{2,1}$ for all
$K \ge 2$.
\end{abstract}

\section{Introduction}

The Brin-Higman-Thompson groups $n G_{k,1}$, for $n \ge 1$ and $k \ge 2$,
are a generalization of the Higman-Thompson group $G_{k,1}$ \cite{Hig74},
and the Brin-Thompson group $n G_{2,1}$ (also called $n V$) 
\cite{Brin1,Brin2}, both of which are generalizations of the Tompson group 
$G_{2,1}$ (also called $V$) \cite{Th, CFP}.
The Brin-Higman-Thompson groups proper, i.e., $n G_{k,1}$ for $n \ge 2$ and 
$k \ge 3$, are the focus of the present paper; in the literature they have
been studied in a number of publications (one of the earliest is 
\cite{DicksMartinez}).

The groups $n G_{2,1}$ and $G_{k,1}$ have many remarkable properties. They 
are finitely presented (for $n G_{2,1}$ this was proved by Brin 
\cite{Brin1, Brin2} and Hennig and Mattucci \cite{HennigMattucci}; for 
$G_{k,1}$ this was proved by Higman \cite{Hig74}). The groups $n G_{2,1}$ 
are simple (Brin \cite{Brin3}), and non-isomorphic for different $n$ (Brin
\cite{Brin1}, Bleak and Lanoue \cite{BleakLanoue}). All $G_{k,1}$ are 
non-isomorphic for different $k$ (Higman \cite{Hig74}, Pardo \cite{Pardo}).
All the groups $G_{k,1}$ are embeddable into one another (\cite{Hig74} and
\cite{MatteBon,BiEmbHTarx,BiEmbHT}); and $(n+1) \, G_{2,1}$ does not 
embed into $n G_{2,1}\,$ (N.\ Matte Bon \cite[Coroll.\ 11.20]{MatteBon}).  
The word problem of $G_{k,1}$ is co-contextfree (Lehnert and Schweitzer
\cite{LehnSchw}), and the word problem of $n G_{2,1}$ is {\sf coNP}-complete
\cite{BinG}.

\medskip

In this paper, we show that $n G_{k,1}$ is finitely generated 
(Theorem \ref{THMnGk1FinGen}), and that the word problem 
over a finite generating set is {\sf coNP}-complete (Theorem \ref{THMwp}).
We show that if $\,K = 1 + (k-1)\,N\,$ for some $N \ge 1$, then
$n G_{K,1} \le n G_{k,1}$; so $n G_{K,1} \le n G_{2,1}$ for all $K \ge 2$
(Theorem \ref{THMembKk}).
We also find a subgroup $2 G_{2,1}^{\rm unif}$ of $2 G_{2,1}$ that is
embeddable in $n G_{k,1}$ for all $n \ge 1$ and $k \ge 3$ (Lemma 
\ref{EmbedG21inGk1}); the group $2 G_{2,1}^{\rm unif}$ contains 
representations of all bijective circuits (Lemma \ref{LEMwpSigmaTauF}). 
The question whether $n G_{2,1}$ is embeddable into $n G_{k,1}$ for 
$k \ge 3$, $n \ge 2$, remains open.

\section{Definitions and background}

In this Section we define the Brin-Higman-Thompson groups $n G_{k,1}$, and
present general background material; we closely follow Section 2 of 
\cite{BinG}. Many results and reasonings in the present paper depend on
\cite{BinG}. 

In this paper, just as in \cite{BinG} and \cite{BiEmbHTarx, BiEmbHT}, 
``function'' means partial function.

An alphabet is any finite set. We will use the alphabets $A_1 = \{a_0\}$
$\subseteq$ $A_2 = \{a_0, a_1\}$  $\subseteq$
$A_k = \{a_0, a_1, \, \ldots, a_{k-1}\}$ with $|A_k| = k$, for any integer 
$k \ge 1$.
For an alphabet $A$, the set of all finite sequences of elements of $A$ 
(called {\em strings}) is denoted by $A^*$. For $m \in {\mathbb N}$, $A^m$ 
is the set of strings of length $m$, and for $x \in A^m$ we write
$|x| = m\,$ ({\it length} of $x$); $A^{\le m}$ is the set of strings of
length $\le m$.
The {\em empty string} is denoted by $\e$; and $|\e| = 0$.
The set of all infinite strings indexed by the ordinal $\omega$ is denoted 
by $A^{\omega}$.
For $x_1, x_2 \in A^*$ the {\it concatenation} is denoted by $x_1 x_2$ or
$x_1 \cdot x_2$.
For $S_1, S_2 \subseteq A^*$, the concatenation is $\, S_1 \cdot S_2$  $=$  
$S_1 \, S_2$  $=$  $\{x_1 \cdot x_2 :$ $x_1 \in S_1$ and $x_2 \in S_2\}$.

For $x, p \in A^*$ we say that $p$ is a {\it prefix} of $x\,$ iff
$\,x = pu$ for some $u \in A^*$; this is denoted by $p \le_{\rm pref} x$.
Two strings $x, y \in A^*$ are called {\it prefix-comparable} (denoted
by $x \,\|_{\rm pref}\, y \,$) iff
$\, x \le_{\rm pref} y \,$ or $ \, y \le_{\rm pref} x$.
A {\it prefix code} is a $\,<_{\rm pref}$-antichain, i.e., a set 
$P \subseteq A^*$ such that for all $p_1, p_2 \in P$: 
$\, p_1 \not<_{\rm pref}\, p_2 \,$.
A {\it right ideal} of $A^*$ is, any $R \subseteq A^*$ such that 
$R = R \cdot A^*$. A subset $C \subseteq R$ generates $R$ as a right ideal 
iff $R = C \cdot A^*$.
Every finitely generated right ideal is generated by a unique finite prefix 
code, and this prefix code is the $\subseteq$-minimum generating set of the 
right ideal. 
A {\it maximal prefix code} is a prefix code $P \subseteq A^*$
that is not a strict subset of any other prefix code of $A^*$.

A {\it right ideal morphism} of $A^*$ is a function $f$: $A^* \to A^*$ such 
that for all $x \in {\rm Dom}(f)$ and all $w \in A^*$:
$ \ f(xw) = f(x) \, w$.  Then ${\rm Dom}(f)$  and ${\rm Im}(f)$ are right
ideals.  The prefix code that generates ${\rm Dom}(f)$ is denoted by 
${\rm domC}(f)$, and is called the {\it domain code} of $f$; the prefix code 
that generates ${\rm Im}(f)$ is denoted by ${\rm imC}(f)$, and is called the 
{\it image code}.

\smallskip

The $n$-fold cartesian product \ {\large \sf X}$_{_{i=1}}^{^n} A^* \, $ will
be denoted by $nA^*$, as in \cite{BinG}; similarly, $nA^{\omega}$ $\,=\,$
 \ {\large \sf X}$_{_{i=1}}^{^n} A^{\omega}$.
Multiplication in $\, nA^*$ happens coordinatewise, i.e., $nA^*$ is the
direct product of $n$ copies of the free monoid $A^*$. For $u \in nA^*$ we
denote the coordinates of $u$ by $u_i \in A^*$, for $1 \le i \le n$; i.e.,
$u = (u_1, \, \ldots, u_n)$.

\smallskip

{\it Geometric interpretation:}  We rename $A_k$ to
$\{0, 1, \,\ldots\, , k-1\}$ and use $A_k$ as the digits for the 
representation of integers and fractions in base-$k$ representation. Then
$\,b_N \ldots b_1 b_0 \centerdot b_{-1} \ldots b_{-S}\,$  (with
$\,b_N,\,\ldots, b_1, b_0, b_{-1},\,\ldots, b_{-S} \in A_k$)
represents the rational number
$\, \sum_{i=0}^N b_i k^i \,+\, \sum_{i=1}^S b_{-i} k^{-i}$.
A rational number $r$ has a finite representation in base $k$ iff
$r = a/k^M$ for some $a \in {\mathbb Z}$ and $M \in {\mathbb N}$; these
rational numbers are called {\em $k$-ary rationals}, or base-$k$ rationals.

We use $x \in A_k^{\,*}$ to represent the semi-open interval
$\,[0.x, \ 0.x + k^{-|x|}[\,$ $\,\subseteq\,$  $[0,1]$ $\,\subseteq\,$
${\mathbb R}$; here, $0.x$ is the $k$-ary rational in $[0,1]$ represented in
fractional base-$k$ representation; the length of this interval is 
$k^{-|x|}$, where $|x|$ is the length of the string $x$. We make an 
exception however: when $0.x + k^{-|x|} = 1$ then the interval represented 
by $x$ is $\,[0.x,\, 1]\,$ (closed interval). 
A set $P \subseteq A_k^{\,*}$ is a maximal prefix code iff the intervals
represented by the strings in $P$ are a partition (tiling) of $[0,1]$.

More generally, $x = (x_1, \ldots, x_n) \in n A_k^{\,*}$ represents
the hyperrectangle
$\,${\large \sf X}$_{_{i=1}}^{^n} [0.x_i, \ 0.x_i + k^{-|x_i|}[ \,$
$\subseteq$  $[0,1]^n$
(except that ``$0.x_i + k^{-|x_i|}[$'' is replaced by ``$1]$'' if
$\,0.x_i + k^{-|x_i|} = 1$).
The measure, in ${\mathbb R}^n$, of this hyperrectangle is
 \ $k^{-(|x_1| \, + \ \ldots \ + \, |x_n|)}$.
In particular, $(\e)^n$ represents $[0,1]^n$ and has measure 1.

This geometric interpretation gives a translation between the description 
the Thompson groups and the Brin-Thompson groups as given in 
\cite{CFP, Brin1, Brin2}, and the string-based description that we use here
and in \cite{BinG}.  We use the string-based approach because it makes it 
easier to study algorithms and computational complexity. 

\medskip

In $nA^*$, the prefix order is generalized to the {\em initial factor 
order}, defined for $u, v \in nA^*$ by $u \le_{\rm init} v$ iff there 
$u x = v$ for some $x \in nA^*$.
Clearly, $u \le_{\rm init} v$ in $nA^*$ iff
$u_i \le_{\rm pref} v_i$ for all $i= 1, \, \ldots, n$.
An {\em initial factor code} is a set $S \subseteq nA^*$ such
that no two different elements of $S$ are $\le_{\rm init}$-comparable (i.e.,
a $\,<_{\rm init}$-antichain).

An important way in which $nA^*$ with $n \ge 2$ differs from $A^*$ concerns 
the {\em join} operation with respect to $\,\le_{\rm init}$.
For all $n$, the join of $u,v \in nA^*$ is defined by $ \ u \vee v$ $\,=\,$
$\min_{\le_{\rm init}}\{z \in nA^* :$
$ \, u \le_{\rm init} z$ and $v \le_{\rm init} z\}\,$; 
$u \vee v$ does not always exist.  A set $S \subseteq nA^*$ is {\em joinless} 
(also called a {\em joinless code}) iff no two elements of $S$ have a join. 
A set $S \subseteq nA^*$ is a {\em maximal} joinless code iff
$\,S$ is $\subseteq$-maximal among the joinless codes of $nA^*$.

\medskip

{\it Geometric interpretation}, continued:
For $u, v \in nA^*$, $v \le_{\rm init} u$ holds iff
the hyperrectangle $u$ is contained in the hyperrectangle $v\,$ (i.e.,
$\le_{\rm init}$ corresponds to $\supseteq$).  The join $u \vee v$
represents the hyperrectangle obtained by intersecting the
hyperrectangles $u$ and $v$ (so $\vee$ corresponds to $\cap$). 
Joinlessness of a set $C \subseteq nA^*$ means that every two hyperrectangles
in $C$ are {\em disjoint}. A set $C \subseteq nA^*$ is a maximal joinless 
code iff the corresponding set of hyperrectangles is a partition (tiling) of 
$[0,1]^n$.
In an initial factor code, $<_{\rm init}$-incomparability means that no
hyperrectangle in the code is strictly contained in another one.

\medskip

\noindent {\bf \cite[Lemma 2.5]{BinG} (the join):}
 \ For all $\,u = (u_1, \, \ldots, u_n), \ v = (v_1, \, \ldots, v_n)$
$\in nA^*$, the join $\,u \vee v \,$ exists iff
for all $i = 1, \, \ldots, n$: \ $u_i \|_{\rm pref} v_i \ $ in $A^*$.
Moreover, if $\,u \vee v\,$ exists, then
$(u \vee v)_i \,=\, u_i\,$ if $v_i \le_{\rm pref} u_i$; and 
$(u \vee v)_i \,=\, v_i\,$ if $u_i \le_{\rm pref} v_i$.

\medskip

\noindent {\bf \cite[Lemma 2.11]{BinG} (one-step restriction):} 
Let $P \subseteq nA^*$ be a finite set. For any
$\, p = (p_1, \, \ldots, p_n) \in P$ and $i \in \{1, \, \ldots, n\}$, let
  \ $P_{p,i}' \ = \ (P \minus \{p\})$  $ \ \cup \ $
$\{(p_1, \, \ldots, p_{i-1}, \, p_i a, \, p_{i+1}, \, \ldots, p_n) \, : \, $
$a \in A\}$.
Then we have: $P$ is a maximal joinless code iff $\, P_{p,i}' \,$ is
a maximal joinless code.

The set $P_{p,i}'$ is called a {\em one-step restriction} of $P$,
and $P$ is called a {\em one-step extension} of $P_{p,i}'$.

\medskip

The concepts of {\it right ideal morphism}, domain code, and image code in 
$nA^*$ are defined in the same way as for $A^*$. 
We only consider domain and image codes that are {\it joinless} (see 
Subsection 2.3 in \cite{BinG}).  We define the monoid 

\smallskip

$n {\cal RI}^{\sf fin}$  $ \ = \ $
$\{f : \, f$ is a right ideal morphism of $nA^*$ such that $f$ is injective,

\hspace{1.06in} and ${\rm domC}(f)$ and ${\rm imC}(f)$
are {\em finite, maximal, joinless} codes\} .

\smallskip

Every right ideal morphism $f \in n {\cal RI}^{\sf fin}$ is uniquely
determined by its restriction to ${\rm domC}(f)$; this is an obvious
consequence of the fact that $f$ is a right-ideal morphism and
${\rm domC}(f)$ is a joinless code.
So $f$ is determined by the finite function
$\, f$: ${\rm domC}(f) \to {\rm imC}(f)$.
A bijection $F$: $P \to Q$ between finite maximal joinless codes
$P, Q \subseteq nA^* \,$ is called a {\em table}.

It follows from \cite[Lemmas 2.21 and 2.29]{BinG} that
$n {\cal RI}^{\sf fin}$ is indeed closed under composition.

Every function $f \in n {\cal RI}^{\sf fin}$ determines a permutation
$f^{(\omega)}$ of $nA^{\omega}$ as follows. Since ${\rm domC}(f)$ is a 
finite maximal joinless code, for any $w \in nA^{\omega}$ there exists a 
unique $p \in {\rm domC}(f)$ such that $w = pu$ for some 
$u \in nA^{\omega}$;
then we define $f^{(\omega)}$ by $ \ f^{(\omega)}(w) = f(p) \ u$.

\smallskip

\noindent {\bf \cite[Def.\ 2.23]{BinG} (end-equivalence):}
Two right ideal morphisms $f, g \in n {\cal RI}^{\sf fin}$ are
{\em end-equivalent} \ iff \ $f$ and $g$ agree on
$\,{\rm Dom}(f) \,\cap\, {\rm Dom}(g)$.
This will be denoted by $f \equiv_{\rm end} g$.

By \cite[Prop.\ 2.18]{BinG}, $\, {\rm Dom}(f) \cap {\rm Dom}(g)$ is
generated by a joinless code, namely
$\, {\rm domC}(f) \vee {\rm domC}(g)$.

\smallskip

\noindent {\bf \cite[Lemma 2.24]{BinG}:}  \ For all $f, g \in$ 
$n {\cal RI}^{\sf fin}$:
  \ $f \equiv_{\rm end} g$ \ iff \ $f^{(\omega)} = g^{(\omega)}$.

\smallskip

\noindent {\bf \cite[Lemma 2.25]{BinG}:}  \ For all $f_1, f_2 \in $
$n {\cal RI}^{\sf fin}$:
 \ \ $(f_2 \circ f_1)^{(\omega)} = f_2^{(\omega)} \circ f_1^{(\omega)}$.
So the relation $\equiv_{\rm end}$ is a {\em congruence} on
$n {\cal RI}^{\sf fin}$.

\smallskip
 
\noindent {\bf \cite[Def.\ 2.28]{BinG} (Brin-Higman-Thompson group
{\boldmath $n G_{k,1}$}):}
Let $A = \{a_0, \, \ldots, a_{k-1}\}$ and $n \ge 1$, $k \ge 2$.
The Brin-Higman-Thompson group $n G_{k,1}$ is
$ \ n \, {\cal RI}^{\sf fin}/\!\! \equiv_{\rm end}$.
 \ Equivalently, $n G_{k,1}$ is the group determined by the action of
$n \, {\cal RI}^{\sf fin}_A$ on $nA_k^{\,\omega}$.

\medskip

In \cite[Def.\ 2.23]{BinG} (quoted above) we used the congruence
$\equiv_{\rm end}$ between right-ideal morphisms. We define a similar
equivalence between right ideals, as follows.

\begin{defn} \label{DEFequivRI}
 \ Let $X, Y \subseteq n A^*$ be finite sets. We say that $X$ and $Y$ are
{\em end-equivalent}, denoted by $X \equiv_{\rm end} Y$, iff
for every right ideal $R \subseteq n A^*$:
 \ \ $R \,\cap\, X \cdot n A^* = \varnothing$  \ $\Leftrightarrow$ i
 \ $R \,\cap\, Y \cdot n A^* = \varnothing$.
\end{defn}
One can prove (see \cite{BiMonVersion}) that for finite sets $X$ and $Y$
the following are equivalent: \\
$X \equiv_{\rm end} Y$;
 \ \ $X \cdot nA^{\omega} = Y \cdot nA^{\omega}$; \ the symmetric difference
$ \ X \cdot n A^* \,\vartriangle\, Y \cdot n A^* \ $ is finite.

\bigskip

An important element of $2 G_{2,1}$ is the {\em shift} $\sigma$, defined by 
$\,{\rm domC}(\sigma) = \{\e\} \x \{a_0,a_1\}$,
$\,{\rm imC}(\sigma) = \{a_0,a_1\} \x \{\e\}$, and
$\,\sigma(\e, b) = (b, \e)$,
for all $b \in \{a_0,a_1\}$.
Hence, $\sigma(x, \, by) = (bx, \, y)$, for all $b \in \{a_0,a_1\}$, and
$\, (x, y) \in 2 A_2^{\,*}$.

An important element of $G_{2,1}$ is the bit-position transposition 
$\tau_{i,i+1}$, defined by ${\rm domC}(\tau_{i,i+1})$ $\,=\,$ 
${\rm imC}(\tau_{i,i+1})$  $\,=\,$   $A_2^{\,i+1}$, and
$ \ \tau_{i,i+1}(x_1\,\ldots\,x_{i-1}\,x_i\,x_{i+1}\,x_{i+2}\,\ldots \ )$
$\,=\,$
$x_1\,\ldots\,x_{i-1}\,x_{i+1}\,x_i\,x_{i+2}\,\ldots \ \ $. 

Another element of $G_{2,1}$ is the Fredkin gate {\small \sf F}, which on 
an input $x_1 x_2 x_3 \in A_2^{\,3}$ is defined by $\,$
{\small \sf F}$(a_0 x_2 x_3) = a_0 x_2 x_3\,$, and $\,$
{\small \sf F}$(a_1 x_2 x_3) = a_1 x_3 x_2\,$
(see \cite[Section 4.5]{BinG}, \cite{Jordan}). 

We embed $\tau_{i,i+1}$ and {\small \sf F} into $2 G_{2,1}$ as 
$\tau_{i,i+1} \x {\mathbb 1}$, and {\small \sf F} $\x$ ${\mathbb 1}\,$ 
(where ${\mathbb 1}$ is the identity on $A_2^{\,*}$).

These examples are also measure-preserving, as defined next.

\begin{defn} \label{DEFmeasurepres} {\bf (measure).}
 \ The measure of $x \in n A^*$ is defined by
$\,\mu(x) \,=\, |A|^{|x_1| \,+\, \ldots \,+\, |x_n|}\,$ (which is the 
measure of the hyperrectangle in $[0,1]^n$ represented by $x$ in the
geometric interpretation).

A right-ideal morphism $f$ of $n A_k^{\,*}$ is called
{\em measure-preserving} \ iff \ for all $x \in {\rm Dom}(f)$:
$\mu(x) = \mu(f(x))$.

An element $g \in n G_{k,1}$ is called measure-preserving \ iff \ $g$ is
represented by some measure-preserving right-ideal morphism of 
$n A_k^{\,*}$.
\end{defn}
Obviously, $\,\mu(x) = \mu(f(x)) \ $ iff 
 \ $\,\sum_{i=1}^n |x_i| \,=\, \sum_{i=1}^n |(f(x))_i|$. 

It is easy to check that if $F$ is a one-step restriction of $f$ and 
$\mu(x) = \mu(f(x))$ for all $x \in {\rm Dom}(f)$, then 
$\mu(x) = \mu(F(x))$ for all $x \in {\rm Dom}(F)$.
Hence, $f$ is measure-preserving iff $f$ is measure-preserving on 
${\rm domC}(f)$.

It is easy to prove that if $g \in n G_{k,1}$ is represented by some
measure-preserving morphism then all morphisms that represent $g$ are
measure-preserving.

The set of measure-preserving elements of $n G_{k,1}$ is a subgroup.

\begin{defn} \label{DEFdictorder} {\bf (dictionary order).}
 \ We consider an alphabet $A_k = \{a_0,a_1, \,\ldots, a_{k-1}\}$ with a
total order $\,a_0 < a_1 < \,\ldots\, < a_{k-1}$.  The 
{\em dictionary order} on $A_k^{\,*}$ is defined as follows. For all 
$\,u, v \in A_k^{\,*}:$ 

\smallskip

$u \le_{\rm dict} v$ \ \ iff

\smallskip

$u \le_{\rm pref} v$, \ \  or

\smallskip

$u \not\le_{\rm pref} v$, and
there exist $p, s, t \in A_k^{\, *}$ and $\alpha, \beta \in A_k$ such
that $u = p \alpha s$, $v = p \beta t$, and $\alpha < \beta$.
\end{defn}
The dictionary order is a well-known total order on $A_k^{\,*}$.

\bigskip

The (maximal) product codes, defined next, are an important special class
of (maximal) joinless codes.

\begin{defn} \label{DEFcoordjoinless} {\bf (product code).} 

\smallskip

\noindent $\bullet$
A {\em product code} in $n A^*$ is a cartesian product
 \ {\large \sf X}$_{n=1}^n P_i \ $ of $n$ prefix codes $P_i \subseteq A^*$.

\smallskip

\noindent $\bullet$
A {\em maximal product code} in $n A^*$ is a maximal joinless code that is a
product code.

\smallskip

\noindent $\bullet$
A right-ideal morphism $f$ of $n A^*$ is a {\em product code morphism} \ iff
 \ ${\rm domC}(f)$ and ${\rm imC}(f)$ are maximal product codes.

\smallskip

\noindent $\bullet$
An element $g \in n G_{k,1}$ is a {\em product code element} iff $g$ can be 
represented by some product code morphism.
\end{defn}
It is easy to prove that every product code is a {\em joinless code}. 

For any set $S \subseteq n A^*$ and $i \in \{1,\ldots,n\}$, let
$\, S_i = \{s_i \in A^*: s \in S\}$; this is the set of $i$th coordinates of
the elements of $S$. It is easy to prove that if $S$ is a product code then
$S_i$ is a prefix code in $A^*$.

The following is also straightforward:

\begin{lem} \label{LEMmaxProdCode}
 \ A product code {\large \sf X}$_{n=1}^n P_i\,$ in $n A^*$ is maximal, 
as a joinless code, iff every $P_i$ is a maximal prefix code in $A^*$.
  \ \ \  \ \ \ $\Box$
\end{lem}
Not every maximal joinless code is a product code; e.g., $\,C =$
$\{(\e, a_0), (a_0,a_1), (a_1,a_1)\}\,$ is a (maximal) joinless code in 
$2 A_2^{\,*}$, but its set of first coordinates $C_1 = \{\e, a_0, a_1\}$ 
is not a prefix code, hence $C$ is not a product code.

Also, not every maximal joinless code is reachable from $\{\e\}^n$ by 
one-step restrictions (as shown by an example of Lawson and Vdovina
\cite[Ex.\ 12.8]{LawsonVdovina}, but the following Lemma shows that we can 
overcome this limitation.

\begin{lem} \label{LEMdomCimCeps}
 \ Every element of $nG_{k,1}$ can be represented by some right-ideal 
morphism $g$ such that ${\rm domC}(g)$ and ${\rm imC}(g)$ are finite 
maximal joinless codes that are reachable from $\{\e\}^n$ by one-step 
restrictions.
\end{lem}
{\sc Proof.} Let $f$ be a right-ideal morphism of $nA_k^{\,*}$ representing
an element of $nG_{k,1}$.
By \cite[Corollary 2.14(0.1, 0.2)]{BinG} one can apply restriction steps to
$f$ in such a way that the resulting right-ideal morphism $h$ has a
${\rm domC}(h)$ that is reachable from $\{\e\}^n$ by one-step restrictions; 
and $h$ represents the same element of $nG_{k,1}$ as $f$. 
If ${\rm imC}(h)$ is also reachable from $\{\e\}^n$ then we pick $h$ for $g$.
Otherwise, we apply more one-step restrictions to $h$ until the resulting
right-ideal morphism  $g$ has an ${\rm imC}(g)$ that is reachable from 
$\{\e\}^n$ by one-step restrictions (which can be done, by
\cite[Corollary 2.14(0.1)(0.2)]{BinG}); and $g$ represents the same element 
of $nG_{k,1}$ as $h$.
In this process, the domain code remains reachable from $\{\e\}^n$,
since ${\rm domC}(g)$ is reached from ${\rm domC}(h)$, which is reached from
$\{\e\}^n$.
Also, every set reachable from $\{\e\}^n$ by one-step restrictions is a 
finite maximal joinless code (by \cite[Lemma 2.11]{BinG}).  
 \ \ \ $\Box$

\bigskip

\noindent {\bf Notation.} The $n$-tuple of copies of $\e$ in $n A^*$ is 
denoted by $(\e)^n$. And the $n$-tuple that contains one string $u$ in 
coordinate $i$, and copies of $\e$ elsewhere, is denoted by 
$\,((\e)^{i-1}, u, (\e)^{n-i})$, or by $\,(\e^{i-1}, u, \e^{n-i})$.
The set whose only element is $(\e)^n$ is denoted by $\{\e\}^n$ or by
$\{(\e)^n\}$. 

Sometimes we denote the finite set of integers $\{1,\ldots,n\}$ by $[1,n]$.

\section{Finite generation of {\boldmath $\, nG_{k,1}$} }

In this section we prove that $nG_{k,1}$ is finitely generated. 
Presumably the methods from \cite{Brin1,HennigMattucci,Quick} could be 
generalized from $nG_{2,1}$ to $nG_{k,1}$; we do indeed use some results 
from \cite{Brin1,HennigMattucci} here.  But we use a slightly different 
method, which also proves that $nG_{k,1}$ is generated by product code 
elements.

The concept of {\em generating set} is well known: A monoid $M$ has $J$
$\subseteq$ $M$ as a generating set iff every element of $M$ is equal to a
product of elements of $J$. We write $M = \langle J \rangle$. 
When $J \subseteq M$ does not generate $M$, it generates a submonoid of $M$, 
denoted by $\langle J \rangle_M$.
A group $G$ has $J$ $\subseteq$ $G$ as a generating set iff every element of
$G$ is equal to a product of elements of $J \cup J^{-1}$ (where
$J^{-1} = \{j^{-1} : j \in J\}$). Again, we write $G = \langle J \rangle$.
So generation in a group is not the same as generation in a monoid
(because of inversion). Since a group is also a monoid, we can have an
ambiguity; by default, we use group generators in a group (unless we
explicitly say ``monoid generators'').
A monoid $M$ (or a group) is {\em finitely generated} iff there exists a
finite subset $J$ $\subseteq$ $M$ that generates $M$; it doesn't matter
whether monoid or group generation is used here (since $J^{-1}$ is finite
when $J$ is finite).

We also consider generation of a subset of a monoid:

\begin{defn} \label{DEFfgInG} {\bf ({\rm  finite generation of a} set).}
 \ Let $M$ be a monoid and let $S \subseteq M$.
The set $S$ is {\em finitely generated in} $M\,$ iff $\,S$ is a subset of a
finitely generated submonoid of $M$ (i.e.,
$S \subseteq \langle J \rangle_{_M}$ for some finite subset $J \subseteq M$).
\end{defn}

\begin{lem} \label{LEMfingenSet}
 \ If a monoid $M$ is generated by $S \subseteq M$, and the set $S$ is
finitely generated in $M$, then $M$ is finitely generated.
\end{lem}
{\sc Proof.}  Suppose $M = \langle S \rangle_{_M}$, and
$S \subseteq \langle J \rangle_{_M}$ for some finite subset $J \subseteq M$.
Now every element of $M$ can be written as a product of elements of $S$,
which can themselves be written as products of elements of $J$.
So $J$ generates $M$.
 \ \ \ $\Box$

\bigskip

We first define a few infinite subsets of $nG_{k,1}$ whose union generates
$nG_{k,1}$. Then we show that these infinite sets are finitely generated in
$nG_{k,1}$. By Lemma \ref{LEMfingenSet} this implies that $nG_{k,1}$ is
finitely generated.

\begin{defn} \label{DEFpositionsets} \hspace{-0.08in}.

\smallskip

\noindent {\small \rm (1.1)} For every $i \in \{1, \ldots, n\}$ we define the
following subgroup of $nG_{k,1}$:

\medskip

 \ \ \  \ \ \  $nG_{k,1}^{(i)} \ = \ $
 $\{{\mathbb 1}\}^{i-1} \times G_{k,1} \times \{{\mathbb 1}\}^{n-i}$.

\bigskip

\noindent {\small \rm (1.2)} For every $i_1, i_2 \in \{1, \ldots, n\}$ such
that $i_1 \ne i_2$ we define the following subset of $nG_{k,1}$:

\medskip

\begin{minipage}{\textwidth}
 \ \ \  \ \ \ $nG_{k,1}^{(i_1\to i_2)} \ = \ \{g \in nG_{k,1} : \ $
there exists a finite maximal prefix code $P \subseteq A_k^{\,*}$ such that


\hspace{2.1in} ${\rm domC}(g)$  $=$
 $\{\e\}^{i_1-1} \times P \times \{\e\}^{n-i_1}$,

\smallskip

\hspace{2.1in} ${\rm imC}(g)$  $=$
 $\{\e\}^{i_2-1} \times P \times \{\e\}^{n-i_2}$,

\smallskip

\hspace{2.1in}
$g\big((\e)^{i_1-1}, \ p,\, (\e)^{n-i_1}\big)$  $=$
$\big((\e)^{i_2 -1}, \ p, (\e)^{n- i_2}\big)$
 \ for all $p \in P \,\}$
\end{minipage}

\bigskip

\noindent {\small \rm (2)} \ Let ${\cal Q} = \{Q^{(m)}: m \in {\mathbb N}\}$
be a set of finite maximal prefix codes in $A_k^{\,*}$.

\smallskip

\noindent {\small \rm (2.1)} For every $i\in \{1, \ldots, n\}$ we define     
the following subset of $nG_{k,1}$:

\medskip

\begin{minipage}{\textwidth}
 \ \ \  \ \ \  $nG_{k,1}^{(i)}({\cal Q}) \ = \ $
$\{g \in nG_{k,1}^{(i)} : \ $ there exists $m \in {\mathbb N}$
such that

\hspace{2.15in} ${\rm domC}(g)$  $=$   ${\rm imC}(g)$  $=$
$\{\e\}^{i-1} \times Q^{(m)} \times \{\e\}^{n-i} \,\}$.
\end{minipage}

\medskip

 \ \ So $g$ is a permutation of ${\rm domC}(g)$, corresponding
to a permutation of $Q^{(m)}$ (for some $m$).

\medskip

\noindent {\small \rm (2.2)} For every $i_1, i_2 \in \{1, \ldots, n\}$ such
that $i_1 \ne i_2$ we define the following subset of $nG_{k,1}$:

\medskip

\begin{minipage}{\textwidth}
 \ \ \  \ \ \ $nG_{k,1}^{(i_1\to i_2)}({\cal Q}) \ = \ $
$\{g \in nG_{k,1}^{(i_1\to i_2)} : \ $ there exists $m \in {\mathbb N}$
such that

\hspace{2.66in} ${\rm domC}(g)$  $=$
 $\{\e\}^{i_1-1} \times Q^{(m)} \times \{\e\}^{n-i_1}$,

\smallskip

\hspace{2.66in} ${\rm imC}(g)$  $=$
 $\{\e\}^{i_2-1} \times Q^{(m)} \times \{\e\}^{n-i_2}\,\}$.
\end{minipage}
\end{defn}
Intuitively, an element of $nG_{k,1}^{(i_1\to i_2)}$ transports
information (namely a string $p \in P$) from coordinate $i_1$ to coordinate
$i_2$, without changing it; whereas an element of $nG_{k,1}^{(i)}$ transforms 
information within the coordinate $i$; in either case the contents of all 
the other coordinates are left unchanged.

\begin{defn} \label{DEFprefcodePm}
 \ We use the following set 
$\,{\cal P} = \{P^{(m)}: m \in {\mathbb N}\}\,$
of finite maximal prefix codes in $A_k^{\,*}$:

\medskip

 \ \ \  \ \ \ $P^{(m)} \,= \ $
 $\bigcup_{i = 0}^m a_{k-1}^{ \ i} \,A_{k-1}$
 \ \ $\cup$ \ \ $\{a_{k-1}^{ \ m+1}\}$

\medskip

\hspace{0.7in} $= \ $
 $\bigcup_{i = 0}^{m-1} a_{k-1}^{ \ i} \,A_{k-1}$
 \ \ $\cup$ \ \ $a_{k-1}^{ \ m} \,A_k$

\medskip

\hspace{0.7in} $= \ $
 $\{a_{k-1}^{ \ i} \, a_j :\, 0 \le i \le m, \ 0 \le j \le k-2\}$
 \ $\cup$ \ $\{a_{k-1}^{ \ m+1}\}\,$.
\end{defn}
In particular, $P^{(0)} = A_k$.
It is useful to also consider the maximal prefix code
$\,P^{(-1)} = \{\e\}$.
Recall that $A_k = A_{k-1} \cup \{a_{k-1}\}\,$ (for $k \ge 2$).

\smallskip

\noindent We have for all $m \ge -1$:
 \ \ $|P^{(m)}| = 1 + (m+1) \, (k - 1)$.

\begin{defn} \label{DEFsigma} {\bf (generalized shifts).}

\noindent For every $m \ge 0$, let $P^{(m)}$ be as in {\rm Def.\ 
\ref{DEFprefcodePm}}. For every $i_1, i_2 \in \{1, \ldots, n\}$ with 
$i_1 \ne i_2$, let

\medskip

 \ \ \  $\sigma_{i_1\to i_2}^{(m)}$  $\,=\,$
$\{\big( ((\e)^{i_1 -1}, \,p, \,(\e)^{n-i_1}),\, $
$((\e)^{i_2 -1}, \,p, \,(\e)^{n-i_2})\big)\,: \ $
$p \in P^{(m)}\}$   
 \ \ \ \ \ $\big(\in\, nG_{k,1}^{(i_1\to i_2)}({\cal P}) \,\big)$,

\medskip

\noindent where
$\,{\rm domC}(\sigma_{i_1 \to i_2}^{(m)})$  $=$
$\{\e\}^{i_1 -1} \x P^{(m)} \x \{\e\}^{n-i_1}$, and
 \ ${\rm imC}(\sigma_{i_1 \to i_2}^{(m)})$  $=$
$\{\e\}^{i_2 -1} \x P^{(m)} \x \{\e\}^{n-i_2}$.
\end{defn}
{\bf Remarks:}

\smallskip

\noindent (1) We see that $\sigma_{2 \to 1}^{(0)} \in 2G_{k,1}$ is the
{\em shift} $\sigma$ (also called the baker's map).

\smallskip

\noindent (2) Composites of $\sigma_{1\to 2}^{(m)}$ reverse the order of the
arguments. E.g.,
 \ $\sigma_{1\to 2}^{(2)} \sigma_{1\to 2}^{(2)}(pq, \e)$ $\,=\,$
$(\e, qp)$.

\begin{lem} \label{LEMnGki1i2sigma}
 \ Let ${\cal P}$ be as in {\rm Def.\ \ref{DEFprefcodePm}}

\smallskip

\noindent {\small \rm (1)} For all $i_1 \ne i_2$ in $\{1, \ldots, n\}$:
 \ \ $nG_{k,1}^{(i_1\to i_2)}({\cal P}) \,=\,$
$\{\sigma_{i_1\to i_2}^{(m)} \,:\, m \ge 0\}$.

\smallskip

\noindent {\small \rm (2)} For all $m \ge 0$, and all two-by-two different 
$i_1, i_2, i_3$:
 \ \ $\sigma_{i_2\to i_3}^{(m)} \circ \sigma_{i_1\to i_2}^{(m)}(.)$
$=$  $\sigma_{i_1\to i_3}^{(m)}(.)$.
\end{lem}
{\sc Proof.}  This follows immediately from Definitions
\ref{DEFpositionsets}(2.2),  \ref{DEFprefcodePm},
and \ref{DEFsigma}.
 \ \ \ $\Box$

\begin{lem} \label{LEMnGkGenby}
 \ Let $\cal P$ be as in {\rm Def.\ \ref{DEFprefcodePm}}.

\smallskip

\noindent {\small \rm (1)} \ The group $nG_{k,1}$ is generated by

\medskip

 \ \ \  \ \ \ $\bigcup_{i=1}^n nG_{k,1}^{(i)}$ \ \ \ $\cup$
 \ \ \ $\bigcup \,\{nG_{k,1}^{(i_1\to i_2)}({\cal P}) \,: $
$\, i_1, i_2 \in \{1, \ldots, n\}, \ i_1 \ne i_2\}\,$.

\bigskip

\noindent {\small \rm (2)} \ The group $nG_{k,1}$ is generated by

\medskip

 \ \ \  \ \ \ $nG_{k,1}^{(1)}({\cal P})$  \ $\cup$
 \ $\bigcup_{i=1}^{n-1} nG_{k,1}^{(i\to i+1)}({\cal P})\,$;

\medskip

 \ \ it is also generated by

\medskip

 \ \ \  \ \ \ $nG_{k,1}^{(1)}({\cal P})$  \ $\cup$
 \ $\bigcup_{i=2}^n nG_{k,1}^{(1\to i)}({\cal P})\,$.

\end{lem}
Part (2) is a stronger version of part (1).

Intuitively the Lemma says that any transformation in $nG_{k,1}$ can be
obtained by composing (repeatedly) in-place transformations at one
coordinate, and transport of information between coordinates.

\medskip

\noindent
{\sc Proof.} (1) Every element of $nG_{k,1}$ is represented by some
right-ideal morphism $g$ such that from ${\rm domC}(g)$, and also from
${\rm imC}(g)$, one can reach $\{\e\}^n$ by one-step extensions
(Lemma \ref{LEMdomCimCeps}).
Let $i_1$ be the coordinate in which the first extension step is applied to
${\rm domC}(g)$ along a sequence of extension steps from ${\rm domC}(g)$ to
$\{\e\}^n$.
Similarly, let $i_2$ be the coordinate in which the first extension step is
applied to ${\rm imC}(g)$ on the way to $\{\e\}^n$.

Let ${\cal P} = \{P^{(m)}: m \ge 0\}$ be as in Def.\ \ref{DEFprefcodePm}.
We choose $m$ so that $|{\rm domC}(g)| = |{\rm imC}(g)| = |P^{(m)}|$;
this is possible by \cite[Corollary 2.14]{BinG}.
Then we can factor $g$ as $\,g(.) = f_2 h f_1(.)$, where

\smallskip

 \ \ \ ${\rm domC}(f_1) = {\rm domC}(g)$,

\smallskip

 \ \ \ ${\rm imC}(f_1) = {\rm domC}(h)$  $ = $
$\{\e\}^{i_1-1} \x P^{(m)} \x \{\e\}^{n-i_1}$,

\smallskip

 \ \ \ ${\rm imC}(h) = {\rm domC}(f_2)$  $=$
$\{\e\}^{i_2-1} \x P^{(m)} \x \{\e\}^{n-i_2}$,

\smallskip

 \ \ \ ${\rm imC}(f_2) = {\rm imC}(g)$.

\smallskip

\noindent Since ${\rm domC}(g)$ and
$\{\e\}^{i_1-1} \x P^{(m)} \x \{\e\}^{n-i_1}$
have entries in coordinate $i_1$ where a one-step extension can be applied,
we define $f_1$ so that it maps the entries of this one-step extension to
each other; except for this, $f_1$ is an arbitrary bijection between domain
codes.
More precisely, there exists $x = (x_1, \ldots,x_n) \in nA_k^{\,*}$ such
that $\,{\rm domC}(g)$ contains

\smallskip

 \ \ \
$\{(x_1, \ldots, x_{i_1 -1},\, x_{i_1} a_j, \,x_{i_1 +1}, \ldots, x_n):$
$j \in [0,k[\, \}$
 \ $=$  \  $x$  $\cdot$
$(\{\e\}^{i_1-1} \x A_k \x \{\e\}^{n-i_1})$,

\smallskip

\noindent which $f_1$ maps bijectively onto
$ \ \{\e\}^{i_1-1} \times a_{k-1}^{\,m}A_k$ $\times$
$\{\e\}^{n-i_1}$    $\,\subseteq$ $P^{(m)}$ \ by

\medskip

  \ \ \
$f_1\big(x \cdot ((\e)^{i_1-1}, \ a_j,\,(\e)^{n-i_1})\big)$
$ \ = \ $
$((\e)^{i_1-1}, \ a_{k-1}^{ \ m}a_j, \,(\e)^{n-i_1})$,

\medskip

\noindent for all $j \in \{0,1,\ldots,k-1\}$. \ Besides this, $f_1$ is
defined to map
 \ ${\rm domC}(g)$  $\minus$    $x$  $\cdot$
$(\{\e\}^{i_1-1} \x A_k \x \{\e\}^{n-i_1})$
 \ bijectively onto       
 \ $(\{\e\}^{i_1-1} \x P^{(m)} \x \{\e\}^{n-i_1})$
$\minus$
$(\{\e\}^{i_1-1}\x a_{k-1}^{\,m}A_k\x \{\e\}^{n-i_1})$
 \ in an arbitrary way.

Similarly, $f_2$ maps the entries of a one-step extension in coordinate
$i_2$ in $\,\{\e\}^{i_2-1} \times a_{k-1}^{\,m}A_k$$\times$
$\{\e\}^{n-i_2}\,$  to a one-step extension in coordinate $i_2$ in
${\rm imC}(g)$; except for this, $f_2$ is an arbitrary bijection.

Finally, $h$ is defined by $h(.) = f_2^{-1} g f_1^{-1}(.)$; so 
${\rm domC}(h)$ and ${\rm imC}(h)$ are as given above.

Now, by applying a one-step extension to $f_1$ and $f_2$ we obtain smaller
functions, hence by induction, $f_1$ and $f_2$ are generated by the set
given in the Lemma.
The base case of the induction happens for the identity map (when $m = -1$);
the identity function can obviously be factored over $nG_{k,1}^{(i)}$ since
${\mathbb 1} \in nG_{k,1}^{(i)}$.

If $i_1 = i_2$, $\,h$ belongs to $nG_{k,1}^{(i_1)}({\cal P})$.

If $i_1 \ne i_2$ then $\,h \,=\, \sigma_{i_1\to i_2}^{(m)}$  $\circ$
$\big(({\mathbb 1})^{i_1-1} \x \pi \x ({\mathbb 1})^{n-i_1}\big)(.)$,
for some permutation $\pi$ of $P^{(m)}$.
And permutations of a finite maximal prefix code in $A_k^{\,*}$ belong to
$G_{k,1}$, so \ $({\mathbb 1})^{i_1-1} \x \pi \x ({\mathbb 1})^{n-i_1}$
$\in$  $nG_{k,1}^{(i_1)}({\cal P})$. So $h$ is generated by
$\,nG_{k,1}^{(i_1)}({\cal P})$ $\,\cup\,$ $nG_{k,1}^{(i_1\to i_2)}$.

\medskip

\noindent (2) Every element of $nG_{k,1}^{(i_1\to i_2)}({\cal P})$ with
$i_1 < i_2$ is generated by
 \ $\bigcup \,\{nG_{k,1}^{(i\to i+1)} : i_1 \le i < i_2\}$.
Indeed, $\sigma_{i_1\to i_2}^{(m)}(.)$ $\,=\,$
$\sigma_{i_2 - 1 \to i_2}^{(m)}$ $\circ$ $\sigma_{i_2 -2\to i_2 -1}^{(m)}$
$\circ$ \ $\ldots$ \ $\circ$
$\sigma_{i_1 +1\to i_1 +2}^{(m)}$ $\circ$ $\sigma_{i_1\to i_1 +1}^{(m)}(.)$.

The case where $i_1 > i_2$ is similar, since 
$nG_{k,1}^{(i_1\to i_2)}({\cal P})$ consists of
the inverses of the elements of $nG_{k,1}^{(i_2\to i_1)}({\cal P})$.
Recall that we use group generators, so inverses of generators are
automatically available.
And $nG_{k,1}^{(n\to 1)}({\cal P})$ is not needed since it consists of the
inverses of the elements of $nG_{k,1}^{(1\to n)}({\cal P})$.
Hence the set $\,\{nG_{k,1}^{(i_1\to i_2)}({\cal P}): i_1 \ne i_2\}\,$ can 
be replaced by $\,\{nG_{k,1}^{(i\to i+1)}({\cal P}) : 1 \le i < n\}$.

And since $\,\sigma_{i\to i+1}^{(m)}(.)$ $=$
$\sigma_{1\to i+1}^{(m)} \circ \sigma_{i\to 1}^{(m)}(.)$,
the set $\,\{nG_{k,1}^{(i\to i+1)}({\cal P}) : 1 \le i < n\}$ can be
replaced by $\,\{nG_{k,1}^{(1\to i)}({\cal P}) : 2 \le i \le n\}$.

In part (1) of the proof, only the subset $\,nG_{k,1}^{(i)}({\cal P})\,$ of
$nG_{k,1}^{(i)}$ is used. Moreover, $nG_{k,1}^{(i)}({\cal P})$ is generated
by $\,nG_{k,1}^{(1)}({\cal P})$ $ \ \cup \ $
$\bigcup_{j=1}^{i-1} nG_{k,1}^{(j\to j+1)}({\cal P})$; indeed,
for any $({\mathbb 1})^{i-1} \x \pi \x ({\mathbb 1})^{n-i}$ $\in$
$nG_{k,1}^{(i)}({\cal P})$, where $\pi$ is a permutation of $P^{(m)}$, we
have
$\,({\mathbb 1})^{i-1} \x \pi \x ({\mathbb 1})^{n-i}$ $\,=\,$
$\sigma_{1\to i}^{(m)}$ $\circ$  $(\pi \x ({\mathbb 1})^{n-1})$  $\circ$
$\sigma_{i\to 1}^{(m)}(.)$; hence, $nG_{k,1}^{(i)}({\cal P})$  $\subseteq$
$nG_{k,1}^{(1\to i)}({\cal P})$ $\cdot$  $nG_{k,1}^{(1)}({\cal P})$ $\cdot$
$nG_{k,1}^{(i\to 1)}({\cal P})$. So every $nG_{k,1}^{(i)}$ (for
$1 \le i \le n$) is generated by $\,nG_{k,1}^{(1)}({\cal P})$ $\,\cup\,$
$\bigcup_{i=2}^n nG_{k,1}^{(1\to i)}({\cal P})$.
 \ \ \  \ \ \  $\Box$

\begin{cor} \label{CORnGk1genProdCodemorph}
 \ $n G_{k,1}$ is generated by product code elements.
\end{cor}
{\sc Proof.} The right-ideal morphisms that represent elements of 
$nG_{k,1}^{(i)}$ or $nG_{k,1}^{(i_1 \to i_2)}({\cal P})$ are 
explicitly defined as product code morphisms.
 \ \ \ $\Box$

\bigskip

\noindent Since $nG_{k,1}^{(i)}$ is isomorphic to $G_{k,1}$, which is
finitely generated, $nG_{k,1}^{(i)}$ is finitely generated, and
$nG_{k,1}^{(i)}({\cal P})$ is a finitely generated subset of $nG_{k,1}$.
By Lemma \ref{LEMnGkGenby}, in order to prove finite generation
of $n G_{k,1}$ it is sufficient to prove finite generation of the set
$nG_{k,1}^{(1\to i)}({\cal P})$ for every $i \in \{2, \ldots, n\}$.

\begin{lem} \label{LEMP2Gk12fingen}
 \ For every $i \in \{1, \ldots, n-1\}$, the set
$\,n G_{k,1}^{(1\to i)}({\cal P})$ is a finitely generated in
$\,nG_{k,1}\,$ (for ${\cal P}$ given in {\rm Def.\ \ref{DEFprefcodePm}}).
\end{lem}
{\sc Proof.} This is proved in the following Lemmas. We will abbreviate
$\sigma_{1\to i}^{(m)}\,$ by $\,\sigma_i^{(m)}$.

\smallskip

\noindent $\bullet$ \ Lemma \ref{LEMsigmaGk} shows that the subset
$n\,G_{k,1}^{(1\to i)}({\cal P})$ is generated by
$\,n\,G_{k,1}^{(1)} \,\cup\, \{\sigma_i^{(m)} : m \ge 0\}\,$ in 
$n\,G_{k,1}$.

\smallskip

\noindent $\bullet$ \ Obviously, $n\,G_{k,1}^{(1)}$ is isomorphic to
$G_{k,1}$, which is finitely presented \cite{Hig74}.

\smallskip

\noindent $\bullet$ \ By Lemma \ref{LEMsigmaFG},
$\{\sigma_i^{(m)} : m \ge 0\}$ is finitely generated in $n\,G_{k,1}$.
 \ \ \ $\Box$

\begin{lem} \label{LEMsigmaGk}
 \ For $2 \le i \le n$ and $\sigma_i^{(m)}$ ($\,= \sigma_{1\to i}^{(m)}$)
as in {\rm Def.\ \ref{DEFsigma}} we have:

\medskip

$n\,G_{k,1}^{(1\to i)}({\cal P})$ is a set generated in
$\,n\,G_{k,1}$ by
 \  $n\,G_{k,1}^{(1)}({\cal P}) \,\cup\, \{\sigma_i^{(m)} : m \ge 0\}\,$.
\end{lem}
{\sc Proof.} This follows immediately from Lemmas \ref{LEMnGki1i2sigma}(1)
and \ref{LEMnGkGenby}(2).
 \ \ \ $\Box$

\begin{lem} \label{LEMsigmaFG}
 \ The set $\,\{\sigma_i^{(m)} :\, m \ge 0, \ 2 \le i \le n\}\,$ is finitely
generated in $\,n\,G_{k,1}$.
\end{lem}
{\sc Proof.} We use Brin's method \cite{Brin1, Brin2} (outlined more briefly
by Burillo and Cleary \cite{BurilloCleary} and used by Hennig and Mattuci
\cite{HennigMattucci}), which starts with a few infinite families of
generators, and then expresses these in terms of a finite subset.
We only need the generators $\,\{C_{m,i} : m \ge 0, \, 2 \le i \le n\}\,$ 
and $\alpha$, defined below; $C_{m,i}$ generalizes Brin's generator 
``$C_m$'' in $2\,G_{2,1}$, and $\alpha$ generalizes Brin's ``$A_0$'' (we 
changed the name of $A_0$ to prevent mix-ups with the names of our 
alphabets). Brin's work only considers the case where $k=2$.
 \ In this proof we abbreviate $(\e)^r$ by $\e^r$.

\medskip

\noindent {\bf (1)} \ For $m = 0$ and $2 \le i \le n$ we define
$\,C_{0,i} = \sigma_i^{(0)}$. 

\smallskip

\noindent For all $m \ge 1$ we define

\bigskip

\begin{minipage}{\textwidth}
 \ \ \ \ \ \ $C_{m,i}$ \ $=$
 \ $\bigcup_{r=0}^{m-1}$
$\{\big((a_{k-1}^{ \ r} a_j, \,\e^{n-1}),$
$\,(a_{k-1}^{ \ r} a_j, \,\e^{n-1})\big) :\, 0 \le j \le k-2\}$

\medskip

\hspace{0.91in} $\cup$
 \ $\{\big((a_{k-1}^{ \ m} a_j, \,\e^{n-1}),\,$
$(a_{k-1}^{ \ m}, \,\e^{i-2}, \,a_j, \,\e^{n-i})\big)$
$\,: \ 0 \le j \le k-1\}\,$,
\end{minipage}

\medskip

 \ \ \  \ \ \ ${\rm domC}(C_{m,i}) \,=\, P^{(m)} \x \{\e\}^{n-1}$,

\medskip

\begin{minipage}{\textwidth}
 \ \ \  \ \ \ ${\rm imC}(C_{m,i})$  $\,=\,$
$\{(a_{k-1}^{ \ r} a_j, \,\e^{n-1}) : $
    $\, 0 \le j \le k-2, \ 0 \le r \le m-1\}$

\smallskip

\hspace{1.25in} $\cup$ \
$\{(a_{k-1}^{ \ m},\,\e^{i-2},\,a_j,\,\e^{n-i})\big)$
$ :\, 0 \le j \le k-1\}$.
\end{minipage}

\bigskip

\noindent In table form, 

\medskip

$C_{1,i}$ \ $=$
 \ \begin{tabular} {|l|l|}
$a_j, \,\e^{n-1}$ & $a_{k-1} a_j, \,\e^{n-1}$ \\ \hline
$a_j, \,\e^{n-1}$ & $a_{k-1},\,\e^{i-2},\,a_j,\,\e^{n-i}$  \\
{\footnotesize ($0 \le j \le k-2$)} & 
{\footnotesize ($0 \le j \le k-1$)}
\end{tabular} , 

\bigskip

$C_{2,i}$ \ $=$
 \ \begin{tabular} {|l|l|l|}
$a_j, \,\e^{n-1}$ & $a_{k-1} a_j, \,\e^{n-1}$ & 
$a_{k-1}^{ \ 2} a_j, \,\e^{n-1}$ \\ \hline
$a_j, \,\e^{n-1}$ & $a_{k-1} a_j, \,\e^{n-1}$ &
$a_{k-1}^{ \ 2},\,\e^{i-2},\,a_j,\,\e^{n-i}$  \\
{\footnotesize ($0 \le j \le k-2$)} & {\footnotesize ($0 \le j \le k-2$)} &
{\footnotesize ($0 \le j \le k-1$)}
\end{tabular} , 

\bigskip

\noindent and for $m \ge 3$,

\medskip

$C_{m,i}$ \ $=$
 \ \begin{tabular} {|l|l|l|l|l|}
$a_j, \,\e^{n-1}$ & $a_{k-1} a_j, \,\e^{n-1}$ & \ \ \ $\ldots$
 \ \ \ & $a_{k-1}^{m-1} a_j, \,\e^{n-1}$ &
$a_{k-1}^{ \ m} a_j, \,\e^{n-1}$ \\ \hline
$a_j, \,\e^{n-1}$ & $a_{k-1} a_j, \,\e^{n-1}$ &
 \ \ \ $\ldots$ \ \ \ & $a_{k-1}^{m-1} a_j,\,\e^{n-1}$ &
$a_{k-1}^{ \ m},\,\e^{i-2},\,a_j,\,\e^{n-i}$  \\
{\footnotesize ($0 \le j \le k-2$)} & {\footnotesize ($0 \le j \le k-2$)} &
 \ \ \ $\ldots$ \ \ \ & {\footnotesize ($0 \le j \le k-2$)} & {\footnotesize
($0 \le j \le k-1$)}
\end{tabular} .

\bigskip

\noindent We have for all $m \ge 0$:

\bigskip

\noindent {\boldmath $(\star)$} \hspace{1.5in}
$\sigma_i^{(m)}(.) \,=\, C_{0,i}\,C_{1,i} \ \ldots \ C_{m,i}(.)\,$.

\bigskip

\noindent Equation $(\star)$ follows by induction from 
$\sigma_i^{(0)} = C_{0,i}\,$ (which holds by the definition of $C_{0,i}$),
and the inductive step

\smallskip

\hspace{1.55in} $\sigma_i^{(m)}(.) \,=\, \sigma_i^{(m-1)} \, C_{m,i}(.)$,

\smallskip

\noindent for all $m \ge 1$.  Let us verify the inductive step.
Looking at ${\rm domC}(C_{m,i})$ ($= {\rm domC}(\sigma_i^{(m)})$), we
consider the inputs
$(a_{k-1}^{ \ r} a_j, \,\e^{n-1})$ with
$0 \le j \le k-2$, $\,0 \le r \le m-1$; and the inputs
$(a_{k-1}^{ \ m} a_j, \,\e^{n-1})$ with $0 \le j \le k-1$.

\smallskip

\noindent For $0 \le j \le k-2$ and $0 \le r \le m-1$:

\smallskip

$(a_{k-1}^{ \ r} a_j, \,\e^{n-1})$
 \ $\stackrel{C_{m,i}}{\longmapsto}$
 \ $(a_{k-1}^{ \ r} a_j, \,\e^{n-1})$
 \ $\stackrel{ \ \sigma_i^{(m-1)}}{\longmapsto}$
 \ $(\e^{i-1}, \,a_{k-1}^{ \ r} a_j, \,\e^{n-i})$
 \ $=$ \ $\sigma_i^{(m)}(a_{k-1}^{ \ r} a_j, \,\e^{n-1})\,$.

\medskip

\noindent For $0 \le j \le k-1$:

\smallskip

$(a_{k-1}^{ \ m} a_j, \,\e^{n-1})$
 \ $\stackrel{C_{m,i}}{\longmapsto}$
 \ $(a_{k-1}^{ \ m}, \,\e^{i-2}, \,a_j, \,\e^{n-i})$
 \ $\stackrel{ \ \sigma_i^{(m-1)}}{\longmapsto}$
 \ $(\e^{i-1}, \,a_{k-1}^{ \ m} a_j, \,\e^{n-i})$
 \ $=$ \ $\sigma_i^{(m)}(a_{k-1}^{ \ m} a_j, \,\e^{n-1})$.

\bigskip

\noindent {\bf (2)} \ We define

\medskip

\begin{minipage}{\textwidth}
 \ \ \  \ \ \ $\alpha$ \ $=$
 \ $\{\big((a_j,\,\e^{n-1}),\,(a_0 a_j,\,\e^{n-1})\big) :$
$ \, 0 \le j \le k-2\}$

\medskip

\hspace{0.75in} $\cup$
 \ \ $\{\big((a_{k-1} a_0, \,\e^{n-1}), \,(a_0 a_{k-1},
\,\e^{n-1})\big)\}$

\medskip

\hspace{0.75in} $\cup$
 \ $\{\big((a_{k-1} a_j, \,\e^{n-1}), \,(a_j, \,\e^{n-1}   )\big) : \,$
            $1 \le j \le k-1\}\,$,
\end{minipage}

\bigskip

 \ \ \  \ \ \ ${\rm domC}(\alpha)$ \ $=$
 \ $(A_{k-1} \ \cup \ a_{k-1} A_k) \times \{\e\}^{n-1}$,

\medskip

 \ \ \  \ \ \ ${\rm imC}(\alpha)$ \ $=$
 \ $(a_0 A_k \ \cup \ (A_k \minus A_1)) \times \{\e\}^{n-1}$;

\medskip

  \ \ \  \ \ \ $\alpha \,\in\, $
$G_{k,1} \times \{{\mathbb 1}\}^{n-1}$ 
$\,=\, n\,G_{k,1}^{(1)}$.

\bigskip

\noindent In table form,

\bigskip

 \ \ \  \ \ \ $\alpha$ \ $=$
 \ \begin{tabular} {|l|l|l|}
$a_j, \,\e^{n-1}$ & $a_{k-1} a_0, \,\e^{n-1}$ &
$a_{k-1} a_j, \,\e^{n-1}$ \\      \hline
$a_0 a_j, \,\e^{n-1}$ & $a_0 a_{k-1}, \,\e^{n-1}$ &
$a_j, \,\e^{n-1}$ \\
{\footnotesize ($0 \le j \le k-2$)} & & {\footnotesize ($1 \le j \le k-1$)}
\end{tabular} .

\bigskip

\medskip

\noindent We have for all $m \ge 1$:

\bigskip

\noindent {\boldmath $(\star\star)$} \hspace{1.5in}
$C_{m+1,i}(.) \,=\, \alpha^{-m}\, C_{1,i} \ \alpha^m(.)\,$.

\bigskip

\noindent Equation $(\star\star)$ follows by induction from

\smallskip

\hspace{1.6in} $C_{m+1,i}(.) \,=\, \alpha^{-1}\, C_{m,i} \ \alpha(.)\,$.

\smallskip

\noindent Let us verify the latter equation for all $m \ge 1$.
Based on the domain codes of $\alpha$ and of $C_{m+1,i}$ we consider the
following inputs: $ \ (a_j, \e^{n-1})$ with $0 \le j \le k-2$;
 \ $(a_{k-1} a_0, \e^{n-1})$;
 \ $(a_{k-1} a_j, \e^{n-1})$ with $1 \le j \le k-2$;
 \ $(a_{k-1}^{ \ r} a_j, \e^{n-1})$ with $1 \le j \le k-2$,
    $ \ 2 \le r \le m$;
 \ $(a_{k-1}^{ \ m+1} a_j, \e^{n-1})$ with $1 \le j \le k-1$.
These inputs exhaust all possibilities since they form a finite maximal
joinless code.

\medskip

\noindent $\bullet$ \ For $0 \le j \le k-2$:
 \ \ $(a_j, \e^{n-1})$
 \ $\stackrel{\alpha}{\longmapsto}$
 \ $(a_0 a_j, \e^{n-1})$
 \ $\stackrel{C_{m,i}}{\longmapsto}$
 \ $(a_0 a_j, \e^{n-1})$
 \ $\stackrel{\alpha^{-1}}{\longmapsto}$
 \ $(a_j, \e^{n-1})$ \ $=$ \ $C_{m+1,i}(a_j, \e^{n-1})$;

in the application of $C_{m,i}$ we used the fact that
$C_{m,i}(a_0, \e^{n-1}) = (a_0, \e^{n-1})$, when $m \ge 1$.

\medskip

\noindent $\bullet$ \ $(a_{k-1} a_0, \e^{n-1})$
 \ $\stackrel{\alpha}{\longmapsto}$
 \ $(a_0 a_{k-1}, \e^{n-1})$
 \ $\stackrel{C_{m,i}}{\longmapsto}$
 \ $(a_0 a_{k-1}, \e^{n-1})$
 \ $\stackrel{\alpha^{-1}}{\longmapsto}$
 \ $(a_{k-1} a_0, \e^{n-1})$  \ $=$
 \ $C_{m+1}(a_{k-1} a_0, \e^{n-1})$;

in the application of $C_{m,i}$ we again used the fact that
$C_{m,i}(a_0, \e^{n-1}) = (a_0, \e^{n-1})$.

\medskip

\noindent $\bullet$ \ For $1 \le j \le k-2$:

 $(a_{k-1} a_j, \e^{n-1})$
 \ $\stackrel{\alpha}{\longmapsto}$
 \ $(a_j, \e^{n-1})$
 \ $\stackrel{C_{m,i}}{\longmapsto}$
 \ $(a_j, \e^{n-1})$
 \ $\stackrel{\alpha^{-1}}{\longmapsto}$
 \ $(a_{k-1} a_j, \e^{n-1})$  \ $=$
 \ $C_{m+1,i}(a_{k-1} a_j, \e^{n-1})$.

\medskip

\noindent $\bullet$ \ For $\,1 \le j \le k-2$, \ $2 \le r \le m$:

 $(a_{k-1}^{ \ r} a_j, \e^{n-1})$
 \ $\stackrel{\alpha}{\longmapsto}$
 \ $(a_{k-1}^{ \ r-1} a_j, \e^{n-1})$
 \ $\stackrel{C_{m,i}}{\longmapsto}$
 \ $(a_{k-1}^{ \ r-1} a_j, \e^{n-1})$
 \ $\stackrel{\alpha^{-1}}{\longmapsto}$
 \ $(a_{k-1}^{ \ r} a_j, \e^{n-1})$ \

\hspace{4.2in} $=$ \ $C_{m+1,i}(a_{k-1}^{ \ r} a_j, \e^{n-1})$.

\medskip

\noindent $\bullet$ \ For $1 \le j \le k-1$:

 $(a_{k-1}^{ \ m+1} a_j, \e^{n-1})$
 \ $\stackrel{\alpha}{\longmapsto}$
 \ $(a_{k-1}^{ \ m} a_j, \e^{n-1})$
 \ $\stackrel{C_{m,i}}{\longmapsto}$
 \ $(a_{k-1}^{ \ m}, \,a_j)$
 \ $\stackrel{\alpha^{-1}}{\longmapsto}$
 \ $(a_{k-1}^{ \ m+1},\,\e^{i-2},\,a_j,\,\e^{n-i})$

\hspace{4.in} $=$ \ $C_{m+1}(a_{k-1}^{ \ m+1} a_j, \e^{n-1})$.

\medskip

\noindent $\Box$

\bigskip

\noindent {\bf Remark:} \ Equation $(\star\star)$ in the proof of Lemma
\ref{LEMsigmaFG} cannot easily be extended to $m=0$; indeed, 
$\,C_{1,i}(.) \,\ne\, \alpha^{-1}\, C_{0,i} \,\alpha(.)$.

\bigskip

Let $N_{k,1}$ be minimum the number of generators of the Higman-Thompson 
group $G_{k,1}$. This number is not known, except for $N_{2,1} = 2\,$ (by 
Thompson \cite{Th}, Mason \cite{DMason}, and Bleak and Quick 
\cite{BleakQuick}).
An upper bound on $N_{k,1}$ can be obtained from the fact that $G_{k,1}$ is
generated by the elements $g \in G_{k,1}$ with
$|{\rm domC}(g)| \le 1 + 3(k-1)\,$ (by Higman \cite[Lemma 4.2]{Hig74}).
Since $nG_{k,1}^{(1)}$ is isomorphic to $G_{k,1}$, the number of generators
of $nG_{k,1}^{(1)}$ is $N_{k,1}$.

\begin{thm} \label{THMnGk1FinGen} {\bf (finite generation of
{\boldmath $nG_{k,1}$}).}

\smallskip

\noindent The Brin-Higman-Thompson group $nG_{k,1}$ is finitely generated.
Moreover:

\smallskip

\noindent {\small \rm 1.} \ Let $\Gamma_{k,1}$ be a finite generating set
of $G_{\!k,1}$, and let
$ \ \Gamma_{k,1}^{(1)} = \Gamma_{k,1} \x \{{\mathbb 1} \}^{n-1}$.
Then $nG_{k,1}$ is generated by
$ \ \Gamma_{k,1}^{(1)} \,\cup\, \{C_{0,i}, \,C_{1,i} \,:\, 2 \le i \le n\}$,
and by $ \ \Gamma_{k,1}^{(1)}$  $\,\cup\,$
$\{\sigma_i^{(0)}, \,\sigma_i^{(1)} \,:\, 2 \le i \le n\}$.

\medskip

\noindent {\rm 2.} \ Let $N_{k,1}$ be the minimum number of generators of
the Higman-Thompson group $G_{k,1}$.  The minimum number of generators of
$nG_{k,1}$ is $ \ \le\, N_{k,1} + 2 \, (n - 1)\,$.
\end{thm}
The elements of $\,\Gamma_{k,1}^{(1)}$  $\,\cup\,$
$\{\sigma_i^{(0)}, \,\sigma_i^{(1)} :\, 2 \le i \le n\}\,$
are product code elements.

\medskip

\noindent 
{\sc Proof.} The Theorem follows from Lemmas
\ref{LEMnGkGenby} and \ref{LEMP2Gk12fingen}.
Since $nG_{k,1}^{(1)}$ is isomorphic to $G_{k,1}$, the number of generators
of $nG_{k,1}^{(1)}$ is $N_{k,1}$.

It follows from equation $(\star\star)$ in the proof of Lemma         
\ref{LEMsigmaFG} that for each $i \in \{2,\ldots,n\}:$ 
 \ $\{C_{0,i},\, C_{1,i},\, \alpha\}\,$ generates $\{C_{m,i} : m \ge 0\}$. 
Since $\alpha$ is in $nG_{k,1}^{(1)}$, $\,\alpha$ is generated by 
$\Gamma_{k,1}^{(1)}$, hence we have: \ $nG_{k,1}\,$ is generated by 
$ \ \Gamma_{k,1}^{(1)} \,\cup\, \{C_{0,i},\,C_{1,i}\,:\,$ 
$2 \le i \le n\}\,$. 

Since only $\{C_{0,i}, C_{1,i}\}$ is needed, and since
$\,C_{0,i} = \sigma_i^{(0)}\,$ and $\,C_{1,i}$ $=$ 
$(\sigma_i^{(0)})^{-1} \, \sigma_i^{(1)}$, it follows that 
 $\,\Gamma_{k,1}^{(1)}$  $\,\cup\,$  
$\{\sigma_i^{(0)}, \,\sigma_i^{(1)} \,:\, 2 \le i \le n\}\,$
generates $nG_{k,1}$.
 \ \ \ $\Box$

\bigskip

\noindent {\bf Comments on the number of generators:} 
For the special case of the Brin-Thompson groups,
i.e.\ $k=2$, it is known that $n G_{2,1}$ is 2-generated for all $n \ge 2\,$ 
(Martyn Quick \cite{Quick} and Collin Bleak \cite[Acknowledgements]{Quick}).
Theorem \ref{THMnGk1FinGen} only implies that the minimum number of 
generators of $\,nG_{2,1}$ is $\,\le 2n\,$.
It was  previously known that $2 G_{2,1}$ has $\le 8$ generators (Brin
\cite[Prop.\ 6.2]{Brin1}), and that $n G_{2,1}$ has $\le 2n+4$ generators
(Hennig and Mattucci \cite[Theorem 25]{HennigMattucci}).

For $k \ge 3$ it remains unknown what the minimum number of generators of 
$\,nG_{k,1}$ is (in particular for $n=1$).

\section{Embedding {\boldmath $nG_{K,1}$} into {\boldmath $nG_{k,1}$ 
         for certain $K > k$} }

The embedding in Theorem \ref{THMembKk} below generalizes a result of 
Higman \cite[Theorem 7.2]{Hig74} from $G_{k,1}$ to $nG_{k,1}$; just as
Higman's proof, it is based on the idea of coding $A_K$ over $A_k^{\,*}$.

\begin{lem} \label{LEMsizemaxjoin}
 \ For every $n \ge 1$ and $k \ge 2$ we have: 

\smallskip

\noindent $\bullet$ \ Every maximal finite joinless code in $n A_k^{\,*}$ 
has cardinality $1 + (k-1) \, N$, for some $N \ge 0$. 

\smallskip

\noindent $\bullet$ \ For every $N \ge 0$ there exists a maximal 
joinless code in $n A_k^{\,*}$ with cardinality $1 + (k-1) \, N$.
In particular, there exists a maximal prefix code in $A_k^{\,*}$ of 
cardinality $1 + (k-1) \, N$.
\end{lem}
{\sc Proof.} See \cite[Coroll.\ 2.14]{BinG}. An example of a maximal prefix 
code of cardinality $1 + (k-1) \, N\,$ is 

\medskip

 \ \ \ $\bigcup_{j=0}^{N-1} A_{k-1} \, a_{k-1}^{ \ j}$  
 \ $\cup$ \ $\{a_{k-1}^{ \ N}\}\,$.
\hspace{0.6in} $\Box$

\bigskip

\noindent
The following is a classical property of prefix codes and is easy to prove.

\begin{lem} \label{LEMprefcodepref}
 \ Let $P \subseteq A_k^{\,*}$ be any finite maximal prefix code. For any 
$s \in P^*$ (concatenations of code words), and any $u \in A_k^{\,*}$: 
 \ if $su \in P^*$ then $u \in P^*$.  \ \ \  \ \ \ $\Box$ 
\end{lem}

\begin{lem} \label{LEMmaxIffEll} 
 \ Let $C \subseteq n A_k^{\,*}$ be any finite joinless code, and let 
$\, \ell \ge {\rm maxlen}(C)$. Then we have:

$C$ is maximal \ iff \ for every $u \in n A_k^{\,\ell}$ there exists 
$z \in C$ such that  \ $z \le_{\rm init} u$.
\end{lem}
{\sc Proof.} By definition, a joinless code $C$ is maximal iff every 
$v \in n A_k^{\,*}$ has a join with some $z \in C$.

\smallskip

\noindent $[\Rightarrow]$ 
If $C$ is maximal then every $u \in n A_k^{\,\ell}$ has a 
join with some $z \in C$. Since, in addition, 
$|u_i| = \ell = {\rm maxlen}(C)$ we have $z_i \le_{\rm pref} u_i$ for all 
$i \in [1,n]$. Hence, $z \le_{\rm init} u$.

\smallskip

\noindent $[\Leftarrow]$ 
Suppose every $u \in n A_k^{\,\ell}$ has a join with some $z \in C$.
Let $v \in n A_k^{\,*}$ by arbitrary. 

Case 1: \ If $|v_i| \ge \ell$ for all $i \in [1,n]$, then 
$v \ge_{\rm init} u$ for some $u \in n A_k^{\,\ell}$, hence by assumption, 
$v \ge_{\rm init} z$ for some $z \in C$.
So $v$ has a join with an element of $C$.

Case 2: \ If $|v_i| < \ell$ for some $i \in [1,n]$, then let 
$w_i = v_i a_0^{\, \ell - |v_i|}$ for every $i$ such that $|v_i| < \ell$;
and let $w_i = v_i$ if $|v_i| \ge \ell$.  Then 
$w = (w_1, \ldots, w_n) \ge_{\rm init} z$ for some $z \in C$ (by Case 1).
Hence there exists $z \in C$ such that $z_i \le_{\rm pref} w_i$ for all $i$.
If $w_i = v_i$ then $z_i \le_{\rm pref} v_i$.
If $w_i = v_i a_0^{\, \ell - |v_i|}$ then 
$z_i \le_{\rm pref} v_i a_0^{\, \ell - |v_i|}$, which implies 
$z_i \ \|_{\rm pref} \ v_i$. Hence by \cite[Lemma 2.5]{BinG}, 
$z \vee v$ exists. 
 \ \ \ $\Box$

The following generalizes a well-known fact about maximal prefix codes.

\bigskip

\noindent {\bf Notation (coding over $A_k$).} \ Let $S \subseteq A_k^{\,*}$ 
be a finite prefix code. For $K > k \ge 2$, let $c: A_K \to S\,$
be any total function. Then for any string
$v = v_1 \,\ldots\, v_m \in A_K^{\,*}\,$ (with
$v_1, \,\ldots\, , v_m \in A_K$), we define
$\,c(v) = c(v_1) \,\ldots\, c(v_m) \in A_k^{\,*}\,$, i.e., the
concatenation of the strings $c(v_j)$ for $j = 1, \ldots, m$.
And for $q = (q_1, \,\ldots\, , q_n) \in n A_K^{\,*}\,$, we define
$c(q) = (c(q_1), \,\ldots\, , c(q_n)) \in n A_k^{\,*} \ $ (an $n$-tuple
of strings). For $Q \subseteq n A_K^{\,*}\,$ we define
$c(Q) = \{c(q) : q \in Q\}$.

\begin{lem} \label{LEMcodingjlcodes} {\bf (coded joinless code).}
Let $\,c: A_K \to P$ be a bijection, where $P \subseteq A_k^{\,*}$ is a 
finite maximal prefix code; and let $Q \subseteq n A_K^{\,*}$ be a finite 
maximal joinless code. Then $c(Q)$ is a finite maximal joinless code in 
$n A_k^{\,*}$.
\end{lem}
{\sc Proof.}  Let us prove that $c(Q)$ is joinless. 
If $c(q), c(r) \in c(Q)$ have a join then $c(q)_i \ \|_{\rm pref} \ c(r)_i$ 
for all $i \in \{1, \ldots, n\}$ (by \cite[Lemma 2.5]{BinG}). By Lemma 
\ref{LEMprefcodepref}, this implies that $q_i \ \|_{\rm pref} \ r_i\,$.
So, non-joinlessness of $c(Q)$ implies non-joinlessness of $Q$; the result
then follows by contraposition.

\smallskip

\noindent Let us prove that $c(Q)$ is maximal. 
By Lemma \ref{LEMmaxIffEll} it is sufficient to prove that for some 
$\ell \ge {\rm maxlen}(c(Q))$ and for all $u = (u_1, \ldots, u_n)$  $\in$  
$n A_k^{\,\ell}$: $\,r \le_{\rm init} u\,$ for some $r \in c(Q)$.
We choose $\ell = \ell_P \, \ell_Q$, where $\ell_P = {\rm maxlen}(P)$, and
$\ell_Q = {\rm maxlen}(Q)$; then $\ell \ge {\rm maxlen}(c(Q))$.

Since $P$ is a maximal prefix code and since $|u_i| = \ell_P \, \ell_Q\,$ 
we have:  
$\, u_i = p_{1,i} \ \ldots \ p_{\ell_Q,i} \ v_i\,$, for some 
$\,p_{1,i}\, , \ \ldots \ , p_{\ell_Q,i} \in P$, and $v_i \in A_k^{\,*}$. 
Since $P = c(A_K)$, this implies: $\, u_i = c(z_i) \ v_i\,$, for some 
$\,z_i \in A_K^{\,*}$ with $|z_i| = \ell_Q\,$.
So, $c(z) \le_{\rm init} u$, where $z = (z_1, \,\ldots\, , z_n)$.

Since $Q$ is a finite maximal joinless code in $n A_K^{\,*}$ and 
$|z_i| = \ell_Q\,$, Lemma \ref{LEMmaxIffEll} implies: 
$\, q \le_{\rm init} z\,$ for some $q \in Q$.
Hence, $c(q) \le_{\rm init} c(z)$.

So, $\,r \,=\, c(q) \,\le_{\rm init}\, c(z) \,\le_{\rm init}\, u$. 
 \ \ \  \ \ \   \ \ \ $\Box$

\begin{thm} \label{THMembKk}
 \ For every $n \ge 2$, every $k \ge 2$, and every $K > k$ we have: 
 \ If $K = 1 + (k-1) \, N\,$ for some $N \ge 1$, then

\hspace{1.9in}  $nG_{K,1} \,\le\, nG_{k,1}$.

\smallskip

\noindent In particular, for all $K > 2$:  

\hspace{1.9in}  $nG_{K,1} \,\le\, nG_{2,1}$.
\end{thm}
{\sc Proof.} Since $K = 1 + (k-1) \, N$, Lemma \ref{LEMsizemaxjoin} implies 
that there exists a maximal finite prefix code $P \subseteq A_k^{\,*}$ with
$|P| = K$, and a bijection $\,c: A_K \to P$.
Now we map $g \in nG_{K,1}$ to $c(g) \in nG_{k,1}$, defined as follows:
$\,{\rm domC}(c(g)) = c({\rm domC}(g))$, 
 \ ${\rm imC}(c(g)) = c({\rm imC}(g))$, and 
$\,c(g): c(x) \mapsto c(g(x))\,$ for all $x \in {\rm domC}(g)$.
We have the commutative diagram 

\medspace

\begin{minipage}{\textwidth}
\hspace{1.05in} $x \ \ \stackrel{g}{\longrightarrow} \ \ g(x)$ 

\medskip

\hspace{0.92in} $c \, \downarrow$ \hspace{0.47in} $\downarrow \, c$

\smallskip

\hspace{0.9in} $c(x) \ \ \stackrel{c(g)\,}{\longrightarrow} \ \ c(g)(c(x))$ 
$=$   $c(g(x))\,$.
\end{minipage}

\medspace

\medspace

By Lemma \ref{LEMcodingjlcodes}, ${\rm domC}(c(g))$ and ${\rm imC}(c(g))$
are maximal joinless codes, so $c(g) \in nG_{k,1}$.
The coding map $c$ is injective, by definition.
To check that $c$ is a homomorphism, let 
$g = \{(x^1, y^1), \,\ldots\, , (x^m,y^m)\}$, and 
$h = \{(y^1, z^1), \,\ldots\, ,$ $(y^m,z^m)\}\,$ be elements of $nG_{K,1}\,$
(where $x^j, y^j, z^j \in n A_K^{\,*}\,$ for $j = 1, \ldots, m$).
By applying restriction steps we can indeed assume that the domain of $h$ is
the image of $g$.
Then $\,hg(.) = \{(x^1, z^1), \,\ldots\, , (x^m,z^m)\}$.
Hence, $\,c(g) = \{(c(x^1), c(y^1)), \,\ldots\, , (c(x^m), c(y^m))\}$, 
$ \ c(h) = \{(c(y^1), c(z^1)),$  $\,\ldots\, ,$  $(c(y^m), c(z^m))\}$, and
$\,c(h) \, c(g)(.)$ $\,=\,$  $\{(c(x^1), c(z^1)), \,\ldots\, ,$  
$(c(x^m), c(z^m))\}$ $\,=\,$   $c(hg)$.
 \ \ \ $\Box$

\section{Embedding the subgroup {\boldmath $2 G_{2,1}^{\rm unif}$ of 
$2 G_{2,1}$ into $nG_{k,1}$} }

The word problem of $2 G_{2,1}$ is {\sf coNP}-complete \cite{BinG},
so if $2 G_{2,1}$ could be embedded into $2 G_{k,1}$ for all $k \ge 3$, 
it would follow (by Lemma \ref{LEMcomplexityRel}(2)) that the word problem 
of $2 G_{k,1}$ is {\sf coNP}-hard. Unfortunately, we do not know whether
$2 G_{2,1}$ is embeddable into $2 G_{k,1}$ for any $k \ge 3$.
 
In \cite[Subsections 4.5 and 4.6]{BinG} and in Lemma \ref{LEMwpSigmaTauF} 
below we describe a finitely generated subgroup 
$\,\langle \sigma$, $\tau_{1,2} \x {\mathbb 1}$, 
{\small \sf F} $\x$ ${\mathbb 1}\rangle\,$ of $\,2 G_{2,1}$, whose
word problem is also {\sf coNP}-complete. 
In Subsection 5.2  we find a subgroup of $2 G_{2,1}$, called 
$2 G_{2,1}^{\rm unif}$, that contains 
$\,\langle \sigma$, $\tau_{1,2} \x {\mathbb 1}$, 
{\small \sf F} $\x$ ${\mathbb 1}\rangle$, and that can be embedded into 
$2 G_{k,1}$ for every $k \ge 3$.
This implies that the word problem of $2 G_{k,1}$ is {\sf coNP}-hard.

Before defining $n G_{2,1}^{\rm unif}$ and embedding it into $n G_{k,1}$ 
we recall the embedding of $G_{2,1}$ into $G_{k,1}$ for all $k \ge 3\,$ 
\cite{MatteBon, BiEmbHT, BiEmbHTarx};
the embedding of $n G_{2,1}^{\rm unif}$ into $n G_{k,1}$ generalizes the 
methods of \cite{BiEmbHT, BiEmbHTarx} to product code morphisms.

\subsection{Review of the embedding {\boldmath $\,G_{2,1} \le G_{k,1}$} }

For any prefix code $P \subseteq A^*$, let $\,{\sf spref}(P)$  $\,=\,$ 
$\{x \in A^* : (\exists p \in P)[\, x <_{\rm pref} p\,]\,\}\,$ be the set 
of {\em strict prefixes} of the elements of $P$.  
If $P \subseteq A^*$ is a finite {\em maximal} prefix code then the 
cardinalities satisfy 
$ \ |P| \,=\, 1 \,+\, |{\sf spref}(P)| \cdot (|A| -1)\,$.

\smallskip

Notation: \ For any $k \ge 3$,
 \ $A_{[2,k[\,}$  $\,=\,$  $\{a_j \in A_k : 2 \le j < k\}$.

\medskip

\noindent {\bf \cite[Lemma 1.4]{BiEmbHTarx}:}
If $P$ is a finite maximal prefix code in $A_2^{\,*}$ and 
$a_j \in A_{[2,k[\,}$,  then
$ \ P \,\cup\, {\sf spref}(P) \cdot a_j\,$ is a finite maximal prefix code 
in $\{a_0,a_1,a_j\}^*$. 
 \ And $ \ P \,\cup\, {\sf spref}(P) \cdot A_{[2,k[}\,$ is a finite maximal
prefix code in $A_k^{\,*}$.

\medskip

\noindent {\bf \cite[Def.\ 2.2]{BiEmbHTarx}:} 
For a finite set $S \subseteq A^*$, a total function 
$g: A^{\omega} \to A^{\omega}\,$ {\em fixes} $\,S  A^{\omega}$ \ iff 
 \ $g(x) = x\,$ for every $\, x \in S A^{\omega}$.
The {\em fixator} (in $G_{2,1}$) of $S A^{\omega}$ is 
 \ ${\rm Fix}(S A^{\omega}) \,=\, \{g \in G_{2,1} :\, $
$(\forall x \in S A^{\omega}) [\,g(x) = x\,]\,\}$.

\medskip

\noindent {\bf \cite[Lemma 2.3]{BiEmbHTarx}:} 
The group $\,{\rm Fix}(a_0 A_2^{\,\omega})\,$ consists of the 
elements of $G_{2,1}$ that have a table of the following form, where 
$\{u_1, \,\ldots, u_{\ell}\}$ and $\{v_1, \, \ldots, v_{\ell}\}$ are
finite maximal prefix codes in $A_2^{\,*}$:

\medskip

\hspace{1.1in} \begin{tabular} {|l|l|l|l|}
$a_0$ & $a_1 u_1$ & $\dots$ & $a_1 u_{\ell}$ \\ \hline
$a_0$ & $a_1 v_1$ & $\dots$ & $a_1 v_{\ell}$
\end{tabular} .

\bigskip

\noindent {\bf \cite[Def.\ 2.6]{BiEmbHTarx}:} We define the 
{\em $*a_j$-successor} for any $\, a_j \in A_{[2,k[\,}$, and any finite 
maximal prefix code $P \subseteq A_2^{\,*}$ with $|P| \ge 2$.  Let 
$(p_1, \, \ldots, p_{\ell})$ be the list of all elements of $P$ in
increasing dictionary order in $A_2^{\,*}$, where $\ell = |P|$;
then $p_1 \in a_0^{\,*}$.
For every $p_s \in P \minus \{p_1\}$, the {\em $*a_j$-successor}
$\, (p_s)'_j \,$ of $p_s$ is the element of
$\,{\sf spref}(P) \, a_j$, defined as follows, assuming
$\, (p_{s+1})'_j\,, \ \ldots \ , (p_{\ell})'_j \,$ have already been chosen:

\smallskip

 \ \ \  \ \ \  \ \ \ $(p_s)'_j \ = \ $
$\min_{\rm dict}\{\, x a_j \in {\sf spref}(P) \, a_j \ :$
$ \ p_s <_{\rm dict} x a_j \ $ {\rm and}
$ \ x a_j \not\in \{(p_{s+1})'_j\,, \ \ldots \ , (p_{\ell})'_j\} \,\}$,

\smallskip

\noindent where $\min_{\rm dict}$ uses the dictionary order in
$\{a_0,a_1,a_j\}^*$. 

In other words, $(p_s)'_j\,$ is the nearest right-neighbor of $p_s$ in
$\, {\sf spref}(P) \, a_j\,$ that has not yet been associated with another 
$p_m$ for $m > s$. 

\medskip

\noindent {\bf \cite[Lemma 2.7]{BiEmbHTarx} and \cite[Def.\ 2.5]{BiEmbHT}:}
There is a simple formula for the $*a_j$-successor. First, we note that
every element of $\,P \minus \{p_1\}\,$ can be written in a unique way as 
$u a_1 a_0^{\,m}$ for some $u \in A_2^{\,*}$ and $m \ge 0$. Then 
the successor formula is: \ \  $(u a_1 a_0^{\,m})'_j \,=\, u a_j\,$.

\medskip

\noindent {\bf \cite[Lemma 2.8]{BiEmbHTarx}:} 
Let $a_j \in A_{[2,k[\,}$, and let $P \subseteq A_2^{\,*}$ be a finite 
maximal prefix code, ordered as
$\,p_1 <_{\rm dict} \, \dots \, <_{\rm dict} p_{\ell}$, where
$\ell = |P| \ge 2$.  Then:

\smallskip

\noindent {\small $\bullet$} \ $\{(p_2)'_j\,, \ \ldots\, , \,(p_{\ell})'_j\}$
$\,=\,$ ${\sf spref}(P) \ a_j\,$.

\smallskip

\noindent {\small $\bullet$} \ Consider the one-step restriction in which 
$P$ is replaced by
$\,P_r = (P \smallsetminus \{p_r\}) \, \cup \, p_r A_2$, for some
$p_r \in P$.  Then with respect to the prefix code $P_r$,
$\, (p_r a_0)'_j$ and $(p_r a_1)'_j$ are uniquely determined by $p_r$
as follows:

 \ if $\,2 \le r \le \ell\,$ then $\,(p_r a_1)'_j =  p_r a_j\,$ and
$\, (p_r a_0)'_j = (p_r)'_j \,$;

 \ if $r = 1$ then $\,p_1 \in P \cap a_0^{\,*}$, hence
$\,p_1 a_0 \in P_1 \cap a_0^{\,*}$; so $(p_1)'_j$ and
$(p_1 a_0)'_j$ do not exist, whereas $(p_1a_1)'_j = p_1a_j$. 

\medskip

\noindent {\bf \cite[Lemma 2.9]{BiEmbHTarx} (embedding of 
{\boldmath $G_{2,1}$ into $G_{k,1}$}):}  \\      
For every $k \ge 3$ there exists a homomorphic embedding
$\,\iota$: $G_{2,1} \,\hookrightarrow\, G_{k,1}$, defined by
\[ g =  \left[ \begin{array}{lll}
p_1 & \ldots & p_{\ell}  \\
q_1 & \ldots & q_{\ell}
\end{array} \right] \ \ \longmapsto
 \ \ \left[ \begin{array}{llll}
a_0 & a_1p_1 & \ldots & a_1p_{\ell}  \\
a_0 & a_1q_1 & \ldots & a_1q_{\ell}
\end{array} \right]
 \ \ \ \hookrightarrow
\hspace{2.7in}  \]
\[ \hspace{0.in} \iota(g) = \left[ \begin{array}{ll lll lll ll ll llll}
a_0\! &|& a_1p_1 \ \ldots \ a_1p_{\ell} &|
  &(a_1p_1)'_2 \ \ldots \ (a_1p_{\ell})'_2
  &| & \ldots \ & \ \ldots &|
  &(a_1p_1)'_{k-1} \ \ldots \ (a_1p_{\ell})'_{k-1}  \\
a_0\! &| &a_1q_1 \ \ldots \ a_1q_{\ell} &|
  &(a_1q_1)'_2 \ \ldots \ (a_1q_{\ell})'_2
  &| & \ldots \ & \ \ldots &|
  &(a_1q_1)'_{k-1} \ \ldots \ (a_1q_{\ell})'_{k-1}
\end{array} \right], \]

\medskip

\noindent
where $\{p_1, \ldots, p_{\ell}\} = {\rm domC}(g)$ and 
$\{q_1, \ldots, q_{\ell}\} = {\rm imC}(g)$ are finite maximal prefix codes 
in $A_2^{\,*}$. 
Equivalently, the table for $\iota(g)$ is
 
\medskip

 \ $\{(a_0, \, a_0)\}$   $ \ \ \cup \ \ $
$\{(a_1p_r, \, a_1q_r) : 1 \le r \le \ell\}$     $ \ \ \cup \ \ $
$\bigcup_{i=2}^{k-1}$
$\big\{\big((a_1p_r)'_i, \ (a_1q_r)'_i\big) : 1 \le r \le \ell \big\}$.

\medskip

\noindent
In Lemma \ref{LEMiotaMax} below we prove that ${\rm domC}(\iota(g))$ and 
${\rm imC}(\iota(g))$ are maximal prefix codes in $A_k^{\,*}$.

It is easy to see that the embedding $\iota$ is an injective function of 
right-ideal morphisms. To show that it is defined on $G_{2,1}$ one proves 
the following commutation relation, for $r = 1, ..., \ell$:

\smallskip

 \ \ \  \ \ \      
$\iota({\sf restr}_{p_r}(g)) = {\sf restr}_{a_1p_r}(\iota(g))\,$,

\smallskip

\noindent where ${\sf restr}_{p_r}(g)$ is the one-step restriction of $g$ 
over $A_2$ at $p_r \in {\rm domC}(g)$, and ${\sf restr}_{a_1p_r}(\iota(g))$
is a one-step restriction over $A_k$.

To show that $\iota(h) \circ \iota(g)$ $=$ $\iota(h \circ g)$, we 
use the above commutation to make the output row of $\iota(g)$ is equal to 
the input row of $\iota(h)$ (by restrictions).

\begin{lem} \label{LEMiotaMax}
 \ If $P$ is a finite maximal prefix code of $A_2^{\,*}$ then the following 
is a finite maximal prefix code of $A_k^{\,*}$: 
 \ \ \ $\{a_0\}$  $\,\cup\,$   $a_1 P$   $\,\cup\,$ 
$\{(a_1 p)'_j :\, p \in P, \ j \in [2,p[ \,\}\,$.
\end{lem}
{\sc Proof.} If $P$ is a maximal prefix code, then $\{a_0\} \cup a_1 P$ is 
also a maximal prefix code. By \cite[Lemma 1.4]{BiEmbHTarx} (quoted above), 
applied to $\{a_0\} \cup a_1 P$,

\smallskip

$\{a_0\}$  $\,\cup\,$  $a_1 P$  $\,\cup\,$ 
${\sf spref}(\{a_0\} \cup a_1 P) \cdot A_{[2,k[\,}\,$ 

\smallskip

\noindent is a maximal prefix code in $n A_k^{\,*}$.
And by \cite[Lemma 2.8]{BiEmbHTarx},  

\smallskip

${\sf spref}(\{a_0\} \cup a_1 P) \cdot A_{[2,k[\,}\,$

\smallskip

$=\,$ 
$\{ (x)'_j : x \in (\{a_0\} \cup a_1 P) \minus a_0^{\,*},$
$ \ j \in [2,p[ \,\}\,$

\smallskip

$=\,$  $\{ (a_1 p)'_j : p \in P, \ j \in [2,k[ \,\}\,$;

\smallskip

\noindent the latter equality holds since 
$\,(\{a_0\} \cup a_1 P) \minus a_0^{\,*} \ = \ a_1 P$.
This implies the Lemma.  
 \ \ \ $\Box$

\subsection{Product codes}

In this Subsection we study {\em product codes} (Def.\ 
\ref{DEFcoordjoinless}), and {\em uniform product codes} (Def.\ 
\ref{DEFunifProdCode}). 
Then we define {\em uniform restrictions}; these preserve uniform 
product codes. 
And we consider right-ideal morphisms that preserve uniform product codes
and uniform restrictions. These morphisms determine a subgroup 
$n G_{2,1}^{\rm unif}$ of $n G_{2,1}$.

\begin{lem} \label{LEMonestepCartProd2}
 \ If $P \subseteq mA^*$ and $Q \subseteq nA^*$ are finite maximal joinless 
codes that can be reached from $\{\e\}^m$, respectively $\{\e\}^n$, by 
one-step restrictions, then $P \times Q \subseteq (m+n)A^*$ is a finite 
maximal joinless code that can be reached from $\{\e\}^{m+n}$ by one-step 
restrictions.
\end{lem}
{\sc Proof.} This follows from the fact that a one-step restriction is
applied to one coordinate, independently of the other coordinates. Let us
denote reachability of a set $Y$ from a set $X$ via a sequence of one-step restrictions by $X \stackrel{*}{\to} Y$. Then
 \ $\{\e\}^{m+n}$
$ \ \stackrel{*}{\to} \ $   $P \times \{\e\}^n$
$ \ \stackrel{*}{\to} \ $   $P \times Q$.
 \ \ \ $\Box$

\begin{lem} \label{LEMcoowProd} \hspace{-.07in}.

\noindent {\small \rm (1)} \ If $S \subseteq n A^*$ is a product code and if
$S$ can be reached from $\{\e\}^n$ by one-step restrictions, then $S$ is a
maximal product code. 

\smallskip

\noindent {\small \rm (2)} \ If $S \subseteq n A^*$ is a maximal product 
code then $S$ can be reached from $\{\e\}^n$ by one-step restrictions. 
\end{lem}
{\sc Proof.} (1) If $S$ can be reached from $\{\e\}^n$ by one-step 
restrictions then $S$ is a maximal joinless code (by 
\cite[Coroll.\ 2.14]{BinG}). And a product code that is a maximal joinless 
code is a maximal product code (by the remarks after Def.\ 
\ref{DEFcoordjoinless}). 

(2) When $n=1$, it is well known that any finite maximal prefix code can 
be obtained from $\{\e\}$ by one-step restrictions.  For $n \ge 2$, the 
result then follows by Lemma \ref{LEMonestepCartProd2}.
 \ \ \ $\Box$

\begin{defn} \label{DEFcoordwRestr} {\bf (coordinatewise restriction at a
string {\boldmath $u$ in coordinate $i$}).}
 \ Let $S \subseteq n A^*$ be a joinless code, let $i \in \{1,\ldots,n\}$, 
and let $u \in A^*$ be a string that occurs at coordinate $i$ in some 
element of $S$.
The {\em coordinatewise restriction of $S$ at $u$ in coordinate $i$} consists 
of applying a one-step restriction at coordinate $i$ to every $s \in S$ 
such that $s_i = u$.  
In other words: \ $S$ is restricted to 
 
\medskip

  \ \ \ $(S \,\minus\, \{t \in S : t_i = u\})$ $ \ \ \cup \ \ $
$\{t \in S : t_i = u\} \cdot (\{\e\}^{i-1} \x A \x \{\e\}^{n-i})\,$.
\end{defn}
Coordinatewise restrictions are closely related to product codes:

\begin{lem} \label{LEMcoordwRestrCode}
 \ Let $Q \subseteq n A^*$ be a joinless code. 

\smallskip

\noindent {\small \rm (1)} \ \ If $Q \subseteq n A^*$ is a maximal product 
code, and $Q'$ is obtained from $Q$ by a {\em coordinatewise} restriction, 
then $Q'$ is a maximal product code. 

\smallskip

\noindent {\small \rm (2)} \ \ $Q$ is a maximal product code \ iff 
 \ $Q$ can be reached from $\{\e\}^n$ by a finite sequence of 
{\em coordinatewise} restrictions.
\end{lem}
{\sc Proof.} (1) Let $Q =$ {\large \sf X}$_{i=1}^n P_i\,$, where every $P_i$ 
is a maximal prefix code, and consider $i \in \{1,\ldots,n\}$ and 
$p \in P_i$. 
When a coordinatewise restriction is applied to $Q$ in coordinate $i$ at 
$p$, $\,Q$ is replaced by 

\medskip

$Q' \ =$  \ {\large \sf X}$_{j=1}^{i-1} \,P_j$  $\times$
$\big( (P_i \minus \{p\}) \ \cup \ p A \big)$  
$\times$ {\large \sf X}$_{j=i+1}^n P_j \,$,
 
\medskip

\noindent which is a maximal product code; indeed,
$(P_i \minus \{p\}) \ \cup \ p A\,$ is a maximal prefix code in $A^*$ if 
$P_i$ is a maximal prefix code.

\smallskip

\noindent (2) $[\Leftarrow]$ Since $\{\e\}^n$ is a maximal product code, 
(1) implies that any set derived from $\{\e\}^n$ by coordinatewise 
restrictions is a maximal product code. 

\smallskip

\noindent $[\Rightarrow]$ Let $Q =$ {\large \sf X}$_{i=1}^n P_i\,$, as in
the proof of (1).  We  will use the fact that every maximal prefix code 
can be derived from $\{\e\}$ by one-step restrictions in $A^*$. Hence, 
$\{\e\}$ can be reached from $P_n$ be a sequence of one-step extensions.
A one-step extension that is applied to $\,pA \subseteq P_n$, 
can also be applied to $Q$ in the form of a coordinatewise extension; now 
$Q$ becomes \ {\large \sf X}$_{i=1}^{n-1} P_i$  $\times$
$((P_n \minus pA) \,\cup\, \{p\})$. 
The result is a smaller maximal product code. 
By induction on the cardinality of the maximal product code, we conclude
that $\{\e\}^n$ is reachable by coordinatewise extensions.
 \ \ \ $\Box$ 

\bigskip

\noindent {\bf Remark.} A single one-step restriction does not necessarily 
preserve product codes. 
For example, if $Q = 2 A_2 = \{(a_0, a_0),$  $(a_0, a_1),$ $(a_1,a_0),$ 
$(a_1,a_1)\}\,$ is restricted at $(a_0, a_0)$ in coordinate 1, one obtains 
$\,C \,=\, \{(a_0 a_0, a_0),$  $(a_0 a_1, a_0),$  $(a_0, a_1),$  
$(a_1,a_0),$  $(a_1,a_1)\}$, which is not a product code since 
$C_1 = \{a_0 a_0,$  $a_0a_1,$  $a_0,$  $a_1\}\,$ is not a prefix code. 

On the other hand, the coordinatewise restriction of $Q$ at $a_0$ in 
coordinate 1 yields the product code 
$\,\{(a_0a_0, a_0),$  $(a_0 a_1, a_0),$  $(a_0a_0,a_1),$ $(a_0a_1,a_1),$  
$(a_1,a_0),$  $(a_1,a_1)\}$  $\,=\,$  
$\{a_0a_0, a_0a_1, a_1\} \x \{a_0,a_1\}$.

\begin{pro} \label{LEMprodCodeNotGrp}
 \ The set of product code elements in $n G_{k,1}$ is a strict subset of 
$n G_{k,1}$, but it is not a group.
\end{pro}
{\sc Proof.} By Lemma \ref{LEMnGkGenby}, $n G_{k,1}$ is generated by 
product code elements.
The Lemma then follows from the fact that there exist elements of 
$n G_{k,1}$ that are not product code elements, as shown in the next Claim. 

\medskip

\noindent {\bf Claim.} {\it The right-ideal morphism $g$ of $2 A_2^{\,*}$,
given by the following table, does {\rm not} represent a product code 
element of $2 G_{2,1}$.} 

\medskip

 \ \ \   \ \ \  $\,g \,=\,$
\begin{tabular} {|l|l|l|}
$(\e, a_0)$ & $(a_0, a_1)$ & $(a_1,a_1)$ \\ \hline
$(a_0,a_1)$ & $(a_1,a_1)$ & $(\e, a_0)$
\end{tabular} .

\medskip

\noindent Proof of the Claim: If there were a product code morphism that is
$\,\equiv_{\rm end}$-equivalent to $g$ then there would be product code 
morphism that is reachable from $g$ by restriction steps (by Lemma
\ref{LEMcoordwRestrCode}(1) and the non-uniform version of Lemma 
\ref{LEMunifprodPres}(3)).
Let us show that no right-ideal morphism $f$ that is reachable from $g$ by 
restriction steps has product codes for both its domain code and image code. 
If $f$ is obtained from $g$ by a finite sequence of restrictions then 

\smallskip

\noindent (Eq 1)  \hspace{1.1in} 
${\rm domC}(f)$  $\,=\,$  
$(\e, a_0) \cdot U$ $\,\cup\,$  $(a_0, a_1) \cdot V$  $\,\cup\,$  
$(a_1,a_1) \cdot W \,$,

\smallskip

\noindent (Eq 2)  \hspace{1.1in}
${\rm imC}(f)$  $\,=\,$  
$(a_0, a_1) \cdot U$ $\,\cup\,$  $(a_1, a_1) \cdot V$  $\,\cup\,$  
$(\e,a_0) \cdot W \,$,

\smallskip

\noindent for some finite maximal joinless codes $U, V, W$  $\subseteq$ 
$2 A_2^{\,*}$; this follows from the definition of $g$ and 
\cite[Lemma 2.11]{BinG} (quoted earlier). 
Geometrically, $\,(\e, a_0) \cdot U$ is a tiling of $(\e, a_0)$,
 \ $(a_0, a_1) \cdot V$ is a tiling of $(a_0, a_1)$, etc. 

Let us assume, for a contradiction, that $f$ is a product code morphism, 
i.e., $\,{\rm domC}(f)$  $\,=\,$ $P_1 \x P_2$ and $\,{\rm imC}(f)$  $\,=\,$  
$Q_1 \x Q_2$ for some finite maximal prefix codes $P_1, P_2, Q_1, Q_2$ 
$\subseteq$  $A_2^{\,*}$. Since $P_i$ and $Q_i$ are maximal prefix codes, 
we have $P_i = a_0 P_{i,0} \,\cup\, a_1 P_{i,1}$ and 
$Q_i = a_0 Q_{i,0} \,\cup\, a_1 Q_{i,1}$, for some finite maximal prefix 
codes $P_{i,j}, Q_{i,j}$ (for $i=1,2$ and $j=1,2$).
Equality (Eq 1) implies: 

$(\e, a_0) \cdot U$  $=$ $(\e, a_0) \cdot (P_1 \x P_{2,0})\,$,

$(a_0, a_1) \cdot V$ $=$ $(a_0, a_1) \cdot (P_{1,0} \x P_{2,1})\,$

$(a_1, a_1) \cdot V$  $=$  $(a_1, a_1) \cdot (P_{1,1} \x P_{2,1})\,$,

\noindent  hence

$U = P_1 \x P_{2,0}\,$, \ \ $V = P_{1,0} \x P_{2,1}\,$,
 \ \ $W = P_{1,1} \x P_{2,1}\,$.

\noindent Similarly, from equality (Eq 2) we obtain

$U = Q_{1,0} \x Q_{2,1}\,$, \ \ $V = Q_{1,1} \x Q_{2,1}\,$,
 \ \ $W = Q_1 \x Q_{2,0}\,$.

\noindent These six equalities imply 

$P_1 = Q_{1,0}$, \ \ $P_{1,0} = Q_{1,1}$, \ \ $P_{1,1} = Q_1$, 
 \ \ $P_{2,0} = P_{2,1} = Q_{2,1} = Q_{2,0}$.

\noindent Now $P_{1,1} = Q_1 = a_0 Q_{1,0} \,\cup\, a_1 Q_{1,1}$ 
$=$ $a_0 P_1 \,\cup\, a_1 P_{1,0}$.
This implies $a_0 P_1 \subseteq P_{1,1}$, 
hence $a_0 (a_0 P_{1,0} \,\cup\, a_1 P_{1,1})$  $\subseteq$  $P_{1,1}$,
which implies \ $a_0 a_1 P_{1,1}$  $\subseteq$  $P_{1,1}$. 
But since $P_{1,1}$ is a non-empty finite set, the latter inclusion is not 
possible, as it would imply $\,(a_0 a_1)^* \, P_{1,1}$  $\subseteq$  
$P_{1,1}$. 
This completes the proof of the Claim.  
%
 \ \ \ $\Box$

\bigskip

By Prop.\ \ref{LEMprodCodeNotGrp}, the set of all product code morphisms 
does not represent a subgroup of $n G_{2,1}$ that can be embedded into 
$n G_{k,1}$, unless all of $n G_{2,1}$ can be embedded into $n G_{k,1}$ 
(which remains an open problem). 
Hence we will now look at a special kind of maximal product codes, and the 
corresponding morphisms:

\begin{defn} \label{DEFunifProdCode} {\bf (uniform product code).}
 \ A {\em uniform product code} in $n A_2^{\,*}$ is any set of the form
{\Large ${\sf X}_{_{i=1}}^{^n}$}$A_2^{ \ m_i}$, 
 \ for $\,m_1,$  $\ldots,$ $m_n$  $\in$  ${\mathbb N}$. 
\end{defn}
Definition \ref{DEFunifProdCode} implies that every uniform product code is 
a finite maximal product code.

\begin{lem} \label{LEMunifprodIntersect}
 \ The join of two uniform product codes is a uniform product code.

Equivalently, the intersection of two right ideals that are generated by 
uniform product codes, is a right ideal generated by a finite uniform
product code.
\end{lem}
{\sc Proof.} Indeed, \ ({\large \sf X}$_{j=1}^n A^{q_j}$)  $\vee$
({\large \sf X}$_{j=1}^n A^{r_j}$)   $\,=\,$
{\large \sf X}$_{j=1}^n A^{\max\{q_j, r_j\}}$. The statement about
intersections follows by \cite[Prop.\ 2.18(3)]{BinG}.
 \ \ \   \ \ \ $\Box$

\begin{defn} \label{DEFunifRestrict} {\bf (uniform restriction).} 
 \ Let $S \subseteq n A^*$ be a set, and let $i \in \{1,\ldots,n\}$.
The {\em uniform restriction of $S$ in coordinate $i$} 
consists of applying a one-step restriction at coordinate $i$ to every 
$s \in S$. 
In other words, $S$ is restricted to 

\hspace{2.2in}  $\,{\sf restr}_i(S)$ $\,=\,$
$S \cdot (\{\e\}^{i-1} \x A \x \{\e\}^{n-i})\,$.
\end{defn}

\begin{lem} \label{LEMunifProdCodeEpsilon}
 \ For any set $Q \subseteq n A^*$ we have:

\smallskip

\noindent {\small \rm (1)} \ If $Q$ is a uniform product code then the 
{\em uniform restriction} of $Q$ at a coordinate is also a uniform 
product code.

\smallskip

\noindent {\small \rm (2)} \ $Q$ is a uniform product code in $n A^*$ \ iff
 \ $Q$ can be reached from $\{\e\}^n$ by a finite sequence of {\em uniform
restrictions}.

\smallskip

\noindent {\small \rm (3)} \ If $Q$ and $Q'$ are uniform product codes in 
$n A^*$ such that $\,Q' \cdot n A^* \subseteq Q \cdot n A^*$, then $Q'$ is 
reachable from $Q$ by uniform restrictions.   
\end{lem}
{\sc Proof.} (1) Let $Q =$ {\large \sf X}$_{j=1}^n A^{m_j}$.
A one-step uniform restriction of $Q$ at coordinate $i$ is 
$\,Q \cdot (\{\e\}^{i-1} \x A \x \{\e\}^{n-i})$ $=$
{\large \sf X}$_{j=1}^{i-1} A^{m_j}i$  $\times$  $A^{m_i +1}$  $\times$
    {\large \sf X}$_{j=i+1}^n A^{m_j}$.
This is a uniform product code.

(2) Hence (by induction on the number of step), if $Q$ is reached from 
$\{\e\}^n$ by a finite sequence of uniform restrictions, then $Q$ is a 
uniform product code in $n A^*$.

Conversely, part $[\Leftarrow]$ of the proof of Lemma 
\ref{LEMcoordwRestrCode}(2) applies here too; this shows that every uniform 
product code can be reached from $\{\e\}^n$ by a finite sequence of uniform
restrictions.

(3) The assumptions imply that $Q$ and $Q'$ have the from
$Q =$ {\large \sf X}$_{j=1}^n A^{m_j}$, respectively
$Q =$ {\large \sf X}$_{j=1}^n A^{r_j}$, with $m_j \le r_j$ for all $j$.
Now by $\,\sum_{j=1}^n (r_j - m_j)\,$ uniform restriction steps one can 
reach $Q'$ from $Q$.
 \ \ \ $\Box$

\begin{defn} \label{DEFprodcodemorph} {\bf (uniform morphism).}
 \ A right-ideal morphism $f$ of $n A^*$ is a {\em uniform morphism} (or a
uniform product code morphism) \ iff \ ${\rm domC}(f)$ and ${\rm imC}(f)$ 
are both uniform product codes.
\end{defn}

\begin{lem} \label{LEMunifMeaspres}
 \ Every uniform morphism of $n A^*$ is measure-preserving.

On the other hand, there exist measure-preserving product code morphisms 
that are not uniform. 
\end{lem}
{\sc Proof.} Let ${\rm domC}(f) =$ {\large \sf X}$_{i=1}^n A^{m_i}\,$, and
${\rm imC}(f) =$ {\large \sf X}$_{i=1}^n A^{r_i}\,$. Since
$|{\rm domC}(f)| = |{\rm imC}(f)|$ we have:
$|A|^{m_1 \,+\, \ldots \,+\, m_n} = |A|^{r_1 \,+\, \ldots \,+\, r_n}$, which
is equivalent to $\,m_1 \,+\,\ldots\,+\, m_n = r_1 \,+\,\ldots\,+\, r_n$.
Hence for all $x \in {\rm domC}(f)$: $\,\mu(x)$  $\,=\,$
$|A|^{-(m_1 \,+\,\ldots\,+\, m_n)}$  $\,=\,$
$|A|^{-(r_1 \,+\,\ldots\,+\, r_n)}$  $\,=\,$  $\mu(g(x))$.

\medskip

\noindent Proof that the converse is false: An example is

\medskip

$g \,=\,$
\begin{tabular} {|l|l|l|l|}
$a_0,\,a_1a_0$ & $a_0,\,a_1a_1$ & $a_1a_0,\,a_0$ & $a_1a_1,\,a_0$ \\ \hline
$a_1a_1,\,a_0$ & $a_1a_0,\,a_0$ & $a_0,\,a_1a_1$ & $a_0,\,a_1a_0$
\end{tabular} 
{\sf id}$_{A_2^{\,2} \,\x\, A_2^{\,2}}$
\begin{tabular} {|l|}
$a_0,a_0$  \\ \hline
$a_0,a_0$ 
\end{tabular} , 

\medskip

\noindent where 
$\,{\rm domC}(g) \,=\, {\rm imC}(g) \,=\, $
$\{a_0,\, a_1a_0,\, a_1a_1\} \times \{a_0,\, a_1a_0,\, a_1a_1\}$. Then $g$ 
is a product code morphism, and it is easy to check that it is 
measure-preserving. 

\medskip

\noindent {\sf Claim.} If $f$ is a right-ideal morphism obtained from $g$ by 
restrictions (where $g$ is given by the table above), then ${\rm domC}(f)$ 
and ${\rm imC}(f)$ are {\em not} both uniform product codes.

\smallskip

\noindent Proof of the Claim:
Let us assume, for a contradiction, that ${\rm domC}(f)$ and ${\rm imC}(f)$ 
are both uniform product codes. Then, just as in the Claim in the proof of 
Prop.\ \ref{LEMprodCodeNotGrp}, there exist finite maximal prefix codes 
$P_i^{(r)}$ in $A_2^{\,*}$, for $i = 1,2$ and $r = 1, \ldots,9$ ($= |g|$), 
such that

\smallskip

\begin{minipage}{\textwidth}
$({\rm domC}(f))_1 \,=\, $ 
$a_0 P_1^{(1)}$  $\,\cup\,$  $a_0 P_1^{(2)}$  $\,\cup\,$ 
$a_1a_0 P_1^{(3)}$  $\,\cup\,$  $a_1a_1 P_1^{(4)}$ 

\smallskip

\hspace{1.in}
$\,\cup\,$  $a_0a_0 P_1^{(5)}$  $\,\cup\,$   $a_0a_1 P_1^{(6)}$  
$\,\cup\,$   $a_1a_0 P_1^{(7)}$  $\,\cup\,$  $a_1a_1 P_1^{(8)}$  
$\,\cup\,$  $a_0 P_1^{(9)}$  
\end{minipage}

\medskip

\hspace{0.8in} 
$\,=\,$ $A_2^{\,m_1}\,$, \ for some $m_1 \ge 1$;

\medskip

\begin{minipage}{\textwidth}
$({\rm imC}(f))_1 \,=\, $
$a_1a_1 P_1^{(1)}$  $\,\cup\,$  $a_1a_0 P_1^{(2)}$  $\,\cup\,$ 
$a_0 P_1^{(3)}$  $\,\cup\,$  $a_0 P_1^{(4)}$  

\smallskip

\hspace{0.9in}
$\,\cup\,$   $a_0a_0 P_1^{(5)}$  $\,\cup\,$  $a_0a_1 P_1^{(6)}$ 
$\,\cup\,$  $a_1a_0 P_1^{(7)}$  $\,\cup\,$ $a_1a_1 P_1^{(8)}$  
$\,\cup\,$  $a_0 P_1^{(9)}$  
\end{minipage}

\medskip

\hspace{0.7in} $\,=\,$ $A_2^{\,n_1}\,$, \ for some $n_1 \ge 1$.

\medskip

\noindent From the above it follows that $a_0 P_1^{(1)}$ $\subseteq$ 
$A_2^{\,m_1}$, and $a_1a_1 P_1^{(1)}$ $\subseteq$ $A_2^{\,n_1}$. 
Hence all strings in $P_1^{(1)}$ have length $m_1 - 1$ and also $n_1 - 2$; 
so $m_1 = n_1 -1$.
Moreover, $a_1a_0 P_1^{(3)}$ $\subseteq$ $A_2^{\,m_1}$,  
and $a_0 P_1^{(3)}$  $\subseteq$ $A_2^{\,n_1}$.
Hence all strings in $P_1^{(3)}$ have length $m_1 - 2$ and also 
$n_1 - 1$; so $m_1 = n_1 + 1$.
But this implies $n_1 -1 = n_1 + 1$, which is false. This completes the
proof of the Claim.
 \ \ \ $\Box$

\medskip

We now extend Def.\ \ref{DEFunifRestrict} from uniform product codes to 
uniform morphisms.

\begin{defn} \label{DEFunifrestrMorph} {\bf (uniform restriction).} 
 \ If $f$ is a uniform morphism of $n A^*$ then the 
{\em uniform restriction} of $f$ in coordinate $i \in \{1,\ldots,n\}$, is 

\medskip

\hspace{0.7in} ${\sf restr}_i(f)$  $\,=\,$
$\{ \big(x \cdot (\e^{i-1}, a, \e^{n-i}), \ $
  $y \cdot (\e^{i-1}, a, \e^{n-i})\big) : \ a \in A, \ (x,y) \in f\}$.

\medskip

\noindent So, $\,{\rm domC}({\sf restr}_i(f))$  $=$
${\sf restr}_i({\rm domC}(f))$, and
$ \ {\rm imC}({\sf restr}_i(f))$  $=$
${\sf restr}_i({\rm imC}(f))$.
\end{defn}

\begin{lem} \label{LEMunifprodPres}
 \ If $f$ is a uniform morphism of $n A^*$, then:

\smallskip

\noindent {\small \rm (1)} \ Any uniform restriction of $f$ is a 
uniform morphism.

\smallskip

\noindent {\small \rm (2)} \ If $Q \subseteq n A^*$ is a uniform product 
code,  then the restrictions $f|_{Q \cdot n A^*}$ and 
$f^{-1}|_{Q \cdot n A^*}$  are uniform morphisms. 

\smallskip

\noindent {\small \rm (3)} \ For any (uniform) morphisms $f$ and $g$ of 
$n A^*$ we have:

\smallskip

$\,f \equiv_{\rm end} g$ \ iff \ $f \cap g$ is reachable by (uniform)
restrictions from $f$ (and similarly, from $g$).
\end{lem}
{\sc Proof.} (1) Let ${\rm domC}(f) =$ {\large \sf X}$_{j=1}^n A^{d_j}$ and
${\rm imC}(f) =$ {\large \sf X}$_{j=1}^n A^{r_j}$.  
Then the uniform restriction $F$ of $f$ at coordinate $i$ satisfies

\smallskip

$\,{\rm domC}(F)$  $\,=\,$ 
${\rm domC}(f) \cdot (\{\e\}^{i-1} \x A \x \{\e\}^{n-i})$ 
$ \ = \ $ 
{\large \sf X}$_{j=1}^{i-1} A^{d_j}$  $\times$  $A^{d_i +1}$  $\times$ 
    {\large \sf X}$_{j=i+1}^n A^{d_j}\,$,

\smallskip

$\,{\rm imC}(F)$   $ \ = \ $
{\large \sf X}$_{j=1}^{i-1} A^{r_j}$  $\times$  $A^{r_i +1}$  $\times$
    {\large \sf X}$_{j=i+1}^n A^{r_j}$. 

\smallskip

\noindent Indeed, for any $x \in {\rm domC}(f)$ and any $a \in A$:
 \ $f(x \cdot (\e^{i-1}, a, \e^{n-i}))$  $=$
$f(x) \cdot (\e^{i-1}, a, \e^{n-i})$.  
So, ${\rm domC}(F)$ and ${\rm imC}(F)$ are uniform product codes.

\smallskip

\noindent (2) By replacing $Q$ by $\,Q \,\vee\, {\rm domC}(f)\,$ (which is a 
uniform product code, by Lemma \ref{LEMunifprodIntersect}), we can assume 
that $Q \subseteq {\rm Dom}(f)$. 
By Lemma \ref{LEMunifProdCodeEpsilon}(3), $Q$ is reachable from 
${\rm domC}(f)$ by uniform restrictions. Now by using part (1) inductively 
we conclude that $f|_Q$ is a uniform morphism.
The same reasoning can be applied to $f^{-1}$, which is also a uniform 
morphism if $f$ is a uniform morphism. 

\smallskip

\noindent (3) $[\Leftarrow]$ A uniform restrictions consists of a finite
sequence of one-step restrictions, and one-step restrictions preserve
$\equiv_{\rm end}$ (by \cite[Prop.\ 2.26]{BinG}). Hence
$f \equiv_{\rm end} f \cap g$, and $g \equiv_{\rm end} f \cap g$; hence
$f \equiv_{\rm end} g$.

\noindent $[\Rightarrow]$ If $f \equiv_{\rm end} g$ then
$f \equiv_{\rm end} f \cap g$; moreover,  moreover, $f \cap g \subseteq f$.
Hence ${\rm domC}(f \cap g) \equiv_{\rm end} {\rm domC}(f)$, so by Lemma
\ref{LEMunifProdCodeEpsilon}(3), ${\rm domC}(f \cap g)$ is obtained from
${\rm domC}(f)$ uniform restrictions steps. These same restrictions steps 
also yield $f \cap g$ from $f$.
 \ \ \ $\Box$

\begin{lem} \label{LEMunifprodPresProdmorph}
 \ If $f$ is a product code morphism of $n A^*$, then
any uniform restriction of $f$ is also a product code morphism.
\end{lem}
{\sc Proof.} Let ${\rm domC}(f) =$ {\large \sf X}$_{j=1}^n X_j$ and
${\rm imC}(f) =$ {\large \sf X}$_{j=1}^n Y_j\,$, where each $X_j$ and $Y_j$
is a finite maximal prefix code in $A^*$.
Then the uniform restriction $F$ of $f$ at coordinate $i$ satisfies

\smallskip

$\,{\rm domC}(F)$  $\,=\,$
${\rm domC}(f) \cdot (\{\e\}^{i-1} \x A \x \{\e\}^{n-i})$
$ \ = \ $
{\large \sf X}$_{j=1}^{i-1} X_j$  $\times$  $X_i A$  $\times$
    {\large \sf X}$_{j=i+1}^n X_j\,$,

\smallskip

$\,{\rm imC}(F)$   $ \ = \ $
{\large \sf X}$_{j=1}^{i-1} Y_j$  $\times$  $Y_i A$  $\times$
    {\large \sf X}$_{j=i+1}^n Y_j$.
\smallskip

\noindent Indeed, for any $x \in {\rm domC}(f)$ and any $a \in A$:
 \ $f(x \cdot (\e^{i-1}, a, \e^{n-i}))$  $=$
$f(x) \cdot (\e^{i-1}, a, \e^{n-i})$.  

Since concatenation of maximal prefix codes is a maximal prefix code, 
$X_i A$ and $Y_i A$ are maximal prefix codes.  So, ${\rm domC}(F)$ and 
${\rm imC}(F)$ are maximal product codes.
 \ \ \ $\Box$

\begin{lem} \label{LEMunifGroup}
 \ The set of uniform morphisms is closed under composition.
\end{lem}
{\sc Proof.} Let $f_2$ and $f_1$ be two uniform morphisms, and let 
$\,C$  $\,=\,$ ${\rm domC}(f_2) \,\vee\, {\rm imC}(f_1)\,$. By Lemma 
\ref{LEMunifprodPres}, $C$ is a uniform product code.
Let $F_2 = f_2|_C$ and let $F_1 = (f_1^{-1}|_C)^{-1}$. 
Then $\,{\rm domC}(F_2) = {\rm imC}(F_1)$  $\,=\,$ $C$.
By Lemma \ref{LEMunifprodPres}(2), $F_1$ and $F_2$ are uniform morphisms.

Moreover, $\,{\rm domC}(F_2 \circ F_1) = {\rm domC}(F_1)$, which is a uniform 
product code; and $\,{\rm imC}(F_2 \circ F_1) = {\rm imC}(F_2)$, which is a
uniform product code. Hence $F_2 \circ F_1$ ($\, = f_2 \circ f_1$) is a 
uniform morphism.
 \ \ \ $\Box$

\medskip

\noindent By Lemma \ref{LEMunifGroup} the following is a group:

\begin{defn} \label{DEFunifsubnG21} {\bf (the subgroup
{\boldmath $n G_{k,1}^{\rm unif}$}).} 
 \ For any $k \ge 2$, the uniform product code subgroup of $n G_{k,1}$ 
consists of the elements of $n G_{k,1}$ that can be represented by uniform 
morphisms.
The uniform product code subgroup of $n G_{k,1}$ is denoted by
 \ $n G_{k,1}^{\rm unif}$.
\end{defn}

\begin{lem} \label{LEMexANdcouterex}
{\bf (examples and counter-examples for uniform morphisms).}

\smallskip

\noindent {\small \rm (1)} The following are examples of elements of 
$n G_{k,1}^{\rm unif}$: 

\smallskip

\noindent For $n=1$, 
$\,G_{k,1}^{\rm unif} \,=\, {\rm lp}G_{k,1}\,$, defined by 
$\,{\rm lp}G_{k,1} \,=\, \{g \in G_{k,1}: |x| = |g(x)|$ {\rm for all}
$x \in {\rm domC}(g)\}$ \ (the group of length-preserving elements of 
$G_{k,1}$; see {\rm \cite{BiFactor}}).  

\smallskip

\noindent For $n=2=k$: \ $\sigma$ (the shift), and all elements of 
$\,{\rm lp}G_{2,1} \x \{{\mathbb 1}\} \ $
(in particular, $(\tau_{2,1}, {\mathbb 1})$, and
$({\small \sf F}, {\mathbb 1})$, where {\small \sf F} is the Fredkin gate).

\medskip

\noindent {\small \rm (2)} The following product code elements do {\em not} 
belong to $2 G_{2,1}^{\rm unif}$:
 \ All the elements of the Thompson group $F_{2,1}$ and of 
$\, F_{2,1} \x \{{\mathbb 1}\}$ (except for the identity elements 
${\mathbb 1}$ and $({\mathbb 1},{\mathbb 1})$). 
\end{lem}
{\sc Proof.} (1) This is straightforward from the definition of these
elements. \ (2) If $f \in F_{2,1}$ then $f$ is not length-preserving (by 
\cite[Lemma 3.1 and Theorem 3.2]{BiFactor}), whereas ${\mathbb 1}$ is 
obviously length-preserving.
Hence $f \x {\mathbb 1}$ is not measure-preserving, which implies (by
Lemma \ref{LEMunifMeaspres}) that 
$f \x {\mathbb 1} \not\in 2 G_{2,1}^{\rm unif}$. 
 \ \ \ $\Box$

\bigskip

\noindent The following Lemma and Remark will play a crucial role in the 
embedding of $2G_{2,1}^{\rm unif}$ into $2G_{k,1}$.

\begin{lem} \label{LEMprodCompletion} {\bf (completing a maximal product 
code in {\boldmath $n A_2^{\,*}$} to a maximal product code in 
{\boldmath $n A_k^{\,*}$}).}
 \ Let $Q \,=\,$ {\large \sf X}$_{i=1}^n Q_i \ $ be a finite maximal 
product code in $n A_2^{\,*}$; let $k \ge 3$. Then

\medskip

\hspace{0.5in}
{\Large \sf X}$_{i=1}^n (Q_i \ \cup \  {\sf spref}(Q_i) \cdot A_{[2,k[\,})$
 
\medskip

\noindent is a finite maximal product code in $n A_k^{\,*}$. 
\end{lem}
{\sc Proof.} This follows immediately from \cite[Lemma 1.4]{BiEmbHTarx} 
and the remarks after Def.\ \ref{DEFcoordjoinless}.
 \ \ \ $\Box$

\bigskip 

\noindent {\bf Remark.} By \cite[Lemma 2.8]{BiEmbHTarx}, for all 
$i \in \{1,\ldots,n\}$ and $a_j \in A_{[2,k[\,}$:
 
\medskip

\hspace{0.5in}  ${\sf spref}(Q_i \minus \{a_0\}^*) \cdot a_j$ 
$ \ = \ $  $\{ (u)'_j :\, u \in Q_i \minus A_1^{\,*}\}$
  \ \ \ \ ($\, \subseteq A_2^{\,*}$).

\subsection{Embedding of {\boldmath $n G_{2,1}^{\rm unif}$} into 
            {\boldmath $n G_{k,1}$} }

The constructions and proofs follow the methods of 
\cite{BiEmbHT,BiEmbHTarx} very closely.

\medskip

Let $Q \subseteq n A_2^{\,*}$ be a finite maximal product code.
Since for the maximal prefix code $Q_i \subseteq A_2^{\,*}$, the strings in 
$Q_i \cap A_1^{\,*}$ have no successors, we first embed $n G_{2,1}$ into
a fixator that moves only $n$-tuples of strings that start with $a_1$, i.e., 
$n$-tuples in $(a_1)^n \cdot n A_2^{\,*}$.
So we generalize \cite[Def.\ 2.2]{BiEmbHTarx} and 
\cite[Lemma 2.3]{BiEmbHTarx}, and we use an embedding
 \ $\,n G_{2,1}$  $\,\hookrightarrow\,$  
${\sf Fix}\big((n A_2 \minus \{a_1\}^n) \cdot n A_2^{\,\omega}\big)\,$,

\begin{defn} \label{DEFFix} {\bf (fixator).} 
 \ Let $Q \subseteq n A_2^*$ be a finite joinless code. 
The {\em fixator} of $\,Q \cdot n A_2^{\,\omega}$ is

\medskip

\hspace{0.3in}  ${\rm Fix}(Q \cdot n A_2^{\,\omega})$ 
$ \ = \ $
$\{\,g \in n G_{2,1} \,: \ $
$(\forall x \in Q \cdot n A_2^{\,\omega}) [\,g(x) = x\,]\,\}$.
\end{defn}
We abbreviate this by $\,{\rm Fix}(Q)$. 
It is straightforward to check that ${\rm Fix}(Q)$ is a group. We do not 
assume that $Q$ is maximal; if $Q$ is maximal then 
${\rm Fix}(Q) = \{{\mathbb 1}\}$. 

In Def.\ \ref{DEFFix} we use the representation of $n G_{2,1}$ by a total 
action on the $n$-dimensional Cantor space $n A_2^{\,\omega}$.  
Alternatively, if $n G_{2,1}$ is defined by partial actions on $nA_2^{\,*}$ 
and the congruence $\equiv_{\rm end}$, then the fixator is replaced by a 
partial fixator (see \cite{BiEmbHT, BiEmbHTarx} for that formulation).

E.g., for the joinless code $\,Q \,=\, n A_2 \minus \{a_1\}^n\,$ 
we obtain the group 
$\,{\sf Fix}\big((n A_2 \minus \{a_1\}^n) \cdot n A_2^{\,\omega}\big)\,$.
In case $n=2$, $ \ Q \,=\, n A_2 \minus \{a_1\}^n$ $\,=\,$ 
$\{(a_0,a_0),$  $(a_0,a_1),$ $(a_1,a_0)\}$.

\begin{lem} \label{FixTable} \hspace{-0.07in}.

\smallskip

\noindent {\small \rm (1)} The group 
$\,{\sf Fix}(n A_2 \minus \{a_1\}^{\,n})\,$
consists of the elements of $n G_{2,1}$ that have a table of the form

\bigskip

\hspace{0.7in}  \begin{tabular} {|l}
$\,$ \\   
$\,$
\end{tabular}
\hspace{-0.25in} {\sf id}$_{n A_2 \minus \{a_1\}^n}$
\begin{tabular} {|l|l|l|}
$(a_1)^n\,x^{(1)}$ & \ $\dots$ \ & $(a_1)^n\,x^{(\ell)}$ \\ \hline
$(a_1)^n\,y^{(1)}$ & \ $\dots$ \ & $(a_1)^n\,y^{(\ell)}$
\end{tabular}

\bigskip

where $\,\{x^{(1)},\,\dots\,, x^{(\ell)}\}\,$ and 
$\,\{y^{(1)},\,\dots\,, y^{(\ell)}\}\,$ are maximal joinless codes in 
$n A_2^{\,*}$, of cardinality $\ell$.

\medskip

\noindent {\small \rm (2)} The subgroup 
$\,{\sf Fix}(n A_2 \minus \{a_1\}^n)\,$ is isomorphic to $n G_{2,1}$.
\end{lem}
{\bf Proof.} (1) The form of the tables follows immediately from the
definition of $\,{\sf Fix}(n A_2 \minus \{a_1\}^n)$.

\smallskip

\noindent (2) We define an isomorphism
$\,\theta: n G_{2,1} \,\to\, {\sf Fix}(n A_2 \minus \{a_1\}^n)\,$ by

\bigskip

 \ \ \ \  \begin{tabular} {|l|l|l|}
$x^{(1)}$ & \ $\dots$ \ & $x^{(\ell)}$ \\ \hline
$y^{(1)}$ & \ $\dots$ \ & $y^{(\ell)}$
\end{tabular} 
 \ \ \ \  $\longmapsto$
 \ \ \ \  \begin{tabular} {|l}
$\,$ \\
$\,$
\end{tabular}
\hspace{-0.25in} {\sf id}$_{n A_2 \minus \{a_1\}^n}$
\begin{tabular} {|l|l|l|}
$(a_1)^n\,x^{(1)}$ & \ $\dots$ \ & $(a_1)^n\,x^{(\ell)}$ \\ \hline
$(a_1)^n\,y^{(1)}$ & \ $\dots$ \ & $(a_1)^n\,y^{(\ell)}$
\end{tabular} .

\bigskip

\noindent This is obviously a bijection, and it is easy to check that
it is a homomorphism.
 \ \ \ $\Box$

\bigskip

\noindent Besides the fixator (or pointwise stabilizer) above, we will use
(setwise) stabilizers.

\begin{defn} \label{DEFstabilizer} {\bf (stabilizer).} 
 \ Let $h$ be a right-ideal morphism $h$ of $n A_k^{\,*}$ that represents
an element of $n G_{k,1}$, and let $R \subseteq n A_k^{\,*}$ be a right 
ideal. We say that $h$ {\em stabilizes} $R$ \ iff \ for every 
$\,x \in {\rm Dom}(h) \,\cap\, R$: $h(x) \in R$, \ and for every 
$\,y \in {\rm Im}(h) \,\cap\, R$: $h^{-1}(y) \in R$.
\end{defn}

Recall the notation $\,A_{[1,k[\,} = \{a_j \in A_k : 1\le j < k\}$.
The set $n A_k \minus n A_{[1,k[\,}$ consists of the $n$-tuples
of letters in which the letter $a_0$ occurs at least once.
When $n=1$, $\,n A_k \minus n A_{[1,k[\,} = \{a_0\}$, so Lemma 
\ref{LEMprodcode2tok}(2) is a generalization of Lemma \ref{LEMiotaMax}.

\begin{lem} \label{LEMprodcode2tok}
 \ Let \ {\large \sf X}$_{i=1}^n P_i\,$ be a finite maximal product code in
$n A_2^{\,*}$, where $P_i \subseteq A_2^{\,*}$ is a finite maximal prefix 
code (for $i = 1,\ldots,n$).  Then:

\medskip

\noindent {\small \rm (1)} \ The following equality and equivalence hold:

\medskip

 \ \ \ \ \ \ {\large \sf X}$_{i=1}^n \big( a_1 P_i$  
$\,\cup\,$ ${\sf spref}(a_1 P_i) \cdot A_{[2,k[\,} \big)\,$ 

\medskip

 \ \ \  \ \ \ $=$ \ {\large \sf X}$_{i=1}^n \big( a_1 P_i$  $\,\cup\,$
  $\{(a_1 p_i)'_j : p_i \in P_i, \ j \in [2,k[\,\} \big)\,$ 

\medskip

 \ \ \  \ \ \ $\,\equiv_{\rm end} \ n \, A_{[1,k[\,}$.

\bigskip

\noindent {\small \rm (2)} \ The following is a finite maximal joinless 
code in $n A_k^{\,*}$:

\medskip

 \ \ \ \ \ \  $(n A_k \minus n A_{[1,k[\,})$
 \ \ $\cup$
 \ \ {\large \sf X}$_{i=1}^n \big( a_1 P_i$  $\,\cup\,$
 ${\sf spref}(a_1 P_i) \cdot A_{[2,k[\,} \big)\,$.  
\end{lem}
{\sc Proof.} (1) In the proof of Lemma \ref{LEMiotaMax} we saw that 
$ \ {\sf spref}(a_1 P_i) \cdot A_{[2,k[\,}$
$\,=\,$ 
$\{(a_1 p_i)'_j : p_i \in P_i, \ j \in [2,k[\,\}$. The first equality 
now follows immediately.

Recall that the relation $\equiv_{\rm end}$ between finite sets was defined 
in Def.\ \ref{DEFequivRI}.
To prove the $\equiv_{\rm end}$-equivalence in the Lemma we observe first 
that $\,\{a_1\} \equiv_{\rm end} a_1 P_i$, since $P_i$ is a finite maximal 
prefix code in $A_2^{\,*}$.  And 
$\,\{a_0\} \cup a_1 P_i \,\cup\, {\sf spref}(\{a_0\}$  $\,\cup\,$ 
$a_1 P_i) \cdot a_j\,$ is a maximal prefix code of $\{a_0, a_1, a_j\}^*$ (by
\cite[Lemma 1.4]{BiEmbHTarx}), whereas $\,\{a_0\} \cup a_1 P_i\,$ is a
maximal prefix code of $\{a_0, a_1\}^*$.
Hence, 

\smallskip

$\{a_j\} \,=\, \{a_0, a_1, a_j\} \minus \{a_0, a_1\}$ 

\smallskip

$\,\equiv_{\rm end}\,$ 
$\,\big(\{a_0\} \cup a_1 P_i \,\cup\, {\sf spref}(\{a_0\}$  $\,\cup\,$
$a_1 P_i) \cdot a_j \big)$  
 \ $\minus$ \ $\big( \{a_0\} \cup a_1 P_i \big)$ 

\smallskip

$\,=\,$  ${\sf spref}(\{a_0\} \cup a_1 P_i) \cdot a_j$
$ \ = \ $  $\{(a_1 p_i)'_j : p_i \in P_i\}$. 

\smallskip

\noindent Hence, 
 \ $a_1 P_i$  $\,\cup\,$
$\{(a_1 p_i)'_j : p_i \in P_i, \ j \in [2,k[\,\}$
$ \ \equiv_{\rm end} \ $
$\{a_1\} \cup A_{[2,k[\,}$  \ $=$ \ $A_{[1,k[\,}$.

\noindent By taking the cartesian product we obtain the 
$\equiv_{\rm end}$-equivalence in the Lemma.

\smallskip

\noindent (2) Obviously, $n A_k$ is a maximal joinless code in 
$n A_k^{\,*}$.  Since {\large \sf X}$_{i=1}^n \big( a_1 P_i$  $\,\cup\,$
 ${\sf spref}(a_1 P_i) \cdot A_{[2,k[\,} \big)\,$
$ \ \equiv_{\rm end} \ $   $n A_{[1,k[\,}$, as we proved in (1),
it follows that 

\smallskip 

$n A_k \,=\, (n A_k \minus n A_{[1,k[\,}) \,\cup\, n A_{[1,k[\,}$
 \ $\equiv_{\rm end}$ 
 \ $(n A_k \minus n A_{[1,k[\,})$
 \ $\cup$
 \ {\large \sf X}$_{i=1}^n \big( a_1 P_i$  $\,\cup\,$
 ${\sf spref}(a_1 P_i) \cdot A_{[2,k[\,} \big)\,$.

\smallskip

\noindent Hence, the latter is a maximal joinless code in $n A_k^{\,*}$, 
since $n A_k$ is maximal.
 \ \ \ $\Box$

\bigskip

\begin{lem} \label{EmbedG21inGk1} {\bf (embedding).}
 \ For every $k \ge 3$ and $n \ge 2$ there exists a homomorphic embedding
$\, \iota: n G_{2,1}^{\rm unif} \,\hookrightarrow\, n G_{k,1}$.
\end{lem}
{\sc Proof.} We first embed $n G_{2,1}$ into 
${\sf Fix}(n A_2 \minus \{a_1\}^n)$, according to Lemma \ref{FixTable}(2).
Form here on we only consider elements of $n G_{2,1}$ whose domain code 
and image code are uniform product codes. 
We define the following injective function $\iota$.  Let
$g \in n G_{2,1}^{\rm unif}$ be given by a uniform morphism, described by a 
table $\,\{ (x^{(r)},\, y^{(r)}) : 1 \le r \le \ell\}\,$
where $\,{\rm domC}(g) \,=\, \{x^{(1)}, \,\ldots\, , x^{(\ell)}\}$
$=$ $X$ $=\,$ {\large \sf X}$_{i=1}^n X_i\,$,
and $\,{\rm imC}(g) \,=\, \{y^{(1)}, \,\ldots\, , y^{(\ell)}\}$ 
$=$ $Y$ $=\,$ {\large \sf X}$_{i=1}^n Y_i\,$.
Then $\iota(g)$ is defined by 

\bigskip

${\rm domC}(\iota(g)) \ = \ (n A_k \minus n A_{[1,k[\,})$
 \ \ $\cup$ 
 \ \ {\large \sf X}$_{i=1}^n \big( a_1 X_i$  $\,\cup\,$
 ${\sf spref}(a_1 X_i) \cdot A_{[2,k[\,}\big)\,$,

\bigskip

${\rm imC}(\iota(g)) \ = \ (n A_k \minus n A_{[1,k[\,})$
 \ \ $\cup$ 
 \ \ {\large \sf X}$_{i=1}^n \big( a_1 Y_i $
   $\,\cup\,$  ${\sf spref}(a_1 Y_i) \cdot A_{[2,k[\,}\big)\,$,

\bigskip

\begin{minipage}{\textwidth}
$\iota(g) \ = \ {\sf id}_{n A_k \minus n A_{[1,k[\,}}$
 
\medskip

\hspace{0.5in} $\cup$
 \ $\{  (u,v) \, : \,$
$ \ u \in {\rm domC}(\iota(g)) \cap n A_{[1,k[\,} \cdot nA_k^{\,*}$, 
$ \ v \in {\rm imC}(\iota(g)) \cap  n A_{[1,k[\,} \cdot nA_k^{\,*}$, \ and

\smallskip

\hspace{0.82in} for all $r \in \{1,\ldots,\ell\}$, 
$i \in \{1,\ldots,n\}$, and $l \in [2,k[\,$, 

\smallskip

\hspace{1.1in} if $\,u_i = a_1 \, x_i^{(r)} \in a_1 X_i \ $ then 
   $\,v_i = a_1 \, y_i^{(r)}\,$;

\smallskip

\hspace{1.1in} if $\,u_i = (a_1 \, x_i^{(r)})'_l$  $\,\in\,$
  ${\sf spref}(a_1 \, X_i) \cdot a_l \ $ then
$\,v_i = (a_1 \, y_i^{(r)})'_l\,\}\,$
\end{minipage}

\medskip

\begin{minipage}{\textwidth}
\hspace{0.3in} $= \ {\sf id}_{n A_k \minus n A_{[1,k[\,}}$

\smallskip

\hspace{0.5in} $\cup$
 \ $\{ (u, v) : \,$ 
$(\exists (x,y) \in g)(\forall i \in [1,n])[\,(u_i, v_i)$ $=$
$(a_1x_i,\, a_1 y_i)$ 

\smallskip

\hspace{3.0in} or $ \ (\exists l \in [2,k[\,)[\,(u_i, v_i)$  $=$
$((a_1x_i)'_l,\, (a_1y_i)'_l)\,] \,]\, \} \,$.
\end{minipage}

\bigskip

\noindent {\bf Remark.} The cardinality of this table of $\iota(g)$ is 
 \ \ $k^n - (k-1)^n \,+\, |g| \cdot (k-1)^n$. 

Indeed, $|n A_k \minus n A_{[1,k[\,}| = k^n - (k-1)^n$.  
And for each entry $(x,y) \in g$, new entries are obtained as follows:
for every coordinate $i \in [1,n]$ there are $1 + (k-2)$ choices for 
$(u_i, v_i)$ in $\,\{(a_1x_i,\, a_1 y_i)\}$  $\,\cup\,$  
$\{((a_1x_i)'_l,\, (a_1y_i)'_l): l \in [2,k[\, \}$.

\bigskip

\noindent By Lemma \ref{LEMprodcode2tok}, the sets given above for
$\,{\rm domC}(\iota(g))\,$ and $\,{\rm imC}(\iota(g))\,$ are maximal
joinless codes in $n A_k^{\,*}$.
In the formula for $\iota(g)$ above we also use the Remark after Lemma 
\ref{LEMprodCompletion}, namely that 
$ \ {\sf spref}(a_1 X_i) \cdot a_j$ $\,=\,$
$\{ (a_1 \, x_i)'_j : x \in X\}$, and similarly for $Y$. 

\smallskip

Stabilization property of $\iota(g)$: 
We observe that the right-ideal morphism $\iota(g)$ of $n A_k^{\,*}$ 
{\em stabilizes} $n A_2^{\,*}\,$ and $n A_{[1,k[\,}$ (as in Def.\ 
\ref{DEFstabilizer}).
Moreover, $\iota(g)$ stabilizes every set of the form 
 \ {\large \sf X}$_{i=1}^n S_i\,$, where the sets 
$\,S_i, \,\ldots\,, S_n\,$ are chosen arbitrarily in $\,\{A_2^{\,*}\}$  
$\,\cup\,$   $\{A_2^{\,*} a_j : a_j \in A_{[2,k[\,}\}$.

\medskip

The function $\iota$ is well-defined and injective, as a function between
tables. Indeed, for all $r \in \{1, \ldots, \ell\}$, $i \in \{1,\ldots,n\}$,
and $j \in [2,k[\,:$ \ \ $(a_1 \, x^{(r)}_i)'_j\,$ determines $x^{(r)}_i$
(by \cite[Lemmas 2,7 and 2.8]{BinG}), which in turn determines $y^{(r)}_i$ 
(via the table for $g$), which determines $(a_1 \, y^{(r)}_i)'_j$.
The embedding of $n G_{2,1}$ into the fixator introduces a letter $a_1$ into
every string to which $(.)'_j$ is applied; so $(.)'_j$ is well-defined
everywhere it is used.

\bigskip

\noindent  For example, for $n=2$, $\,\iota(g)$ can be described more 
graphically:

\bigskip

\begin{minipage}{\textwidth}
$g \ = \ $
\begin{tabular} {|c|c|c|}
 \ $\dots$ \ & $(x^{(r)}_1, x^{(r)}_2)$ & \ $\dots$ \\ \hline
 \ $\dots$ \ & $(y^{(r)}_1, y^{(r)}_2)$ & \ $\dots$ 
\end{tabular}
 \ \ \ \  $\longmapsto$
 \ \ \ \  \begin{tabular} {|l}
$\,$ \\
$\,$
\end{tabular}
\hspace{-0.25in} {\sf id}$_{2 A_2 \minus \{(a_1,a_1)\}}$
\begin{tabular} {|c|c|c|}
$\dots$ \ & $(a_1\,x^{(r)}_1,\, a_1\,x^{(r)}_2)$ & \ $\dots$ \\ \hline
$\dots$ \ & $(a_1\,y^{(r)}_1,\, a_1\,y^{(r)}_2)$ & \ $\dots$ 
\end{tabular} 
 \ \ \ \  $\longmapsto$

\bigskip

$\iota(g)$ \ $=$
 \ \begin{tabular} {|l}
$\,$ \\
$\,$
\end{tabular}
\hspace{-0.25in} {\sf id}$_{2 A_k \minus 2 A_{[2,k[\,}}$
\begin{tabular} {|c|c|c|c|c|c|}
 \ $\dots$ \ & $(a_1\,x^{(r)}_1,\, a_1\,x^{(r)}_2)$ & \ $\dots$ \ & 
 \ $\dots$ \ & $(a_1\,x^{(r)}_1,\, (a_1\,x^{(r)}_2)'_j)$ & \ $\dots$ 
 \\ \hline
 \ $\dots$ \ & $(a_1\,y^{(r)}_1,\, a_1\,y^{(r)}_2)$ & \ $\dots$ \ &
 \ $\dots$ \ & $(a_1\,y^{(r)}_1,\, (a_1\,y^{(r)}_2)'_j)$ & \ $\dots$
\end{tabular} \ \ \ $\dots$
 
\bigskip

\hspace{1.1in}
$\dots$ \ \ \  \begin{tabular} {|c|c|c|c|c|c|}
 \ $\dots$ \ & $((a_1\,x^{(r)}_1)'_j,\, a_1\,x^{(r)}_2)$ & \ $\dots$ \ &  
 \ $\dots$ \ & $((a_1\,x^{(r)}_1)'_j,\, (a_1\,x^{(r)}_2)'_j)$ & \ $\dots$   
 \\ \hline
 \ $\dots$ \ & $((a_1\,y^{(r)}_1)'_j,\, a_1\,y^{(r)}_2)$ & \ $\dots$ \ &
 \ $\dots$ \ & $((a_1\,y^{(r)}_1)'_j,\, (a_1\,y^{(r)}_2)'_j)$ & \ $\dots$    
\end{tabular} ,
\end{minipage}

\bigskip

\noindent where $1 \le r \le \ell$ and $j \in [2,k[\,$. 

\bigskip

To show that $\iota$ is not just a function from right-ideal morphisms of
$n A_2^{\,*}$ to right-ideal morphisms of $n A_k^{\,*}$, but also a function 
from $n G_{2,1}^{\rm unif}$ into $n G_{k,1}$, we will show that for all 
$i \in \{1,\ldots,n\}$ and all right-ideal morphisms $g$ of $n A_2^{\,*}$:  
The operation of {\em uniform restriction} ${\sf restr}_i(.)$ {\em commutes} 
with $\iota$, as follows.

\bigskip

\noindent $(\star)$ \hspace{1.1in}
$\iota({\sf restr}^{A_2}_i(g))$  $\,=\,$
    ${\sf restr}^{A_2 A_k}_i(\iota(g))\,$,

\bigskip

\noindent where ${\sf restr}^{A_2}_i(g)$ is the uniform restriction, in
coordinate $i$, of the uniform morphism $g$ of $n A_2^{\,*}\,$ (Def.\
\ref{DEFunifrestrMorph}). 
And for a uniform morphism $h$ of $n A_k^{\,*}$,  
$\,{\sf restr}^{A_2 A_k}_i(h)\,$ is defined to be the restriction in 
coordinate $i$ that is uniformly applied to all $(x,y) \in h$ such that 
$x, y \in n A_{[1,k[\,} \cdot n A_k^{\,*}$, and
$x_i, y_i \in A_2^{\,*}$. Note that if $h = \iota(g)$ then either both 
$x_i, y_i \in A_2^{\,*}\,$ or both $x_i, y_i \not\in A_2^{\,*}\,$. More 
precisely, ${\sf restr}^{A_2 A_k}_i(.)$ is defined as follows, where
it is assumed that $h$ has the stabilization properties of $\iota(g)$, 
observed earlier (in particular, $h$ stabilizes $n A_2^{\,*}$ and 
$n A_{[1,k[\,}$):

\medskip

\begin{minipage}{\textwidth}
${\sf restr}^{A_2 A_k}_i(h)$  $\,=\,$
$\{(x,y) \in h : \ x_i, y_i \not\in A_2^{\,*} \ $ or 
$ \ x = y \in (n A_k \minus n A_{[1,k[\,}) \cdot n A_k^{\,*} \,\}$

\medskip

\hspace{1.05in} $\,\cup\,$
$\{\big(x \cdot (\e^{i-1}, a, \e^{n-i}),\, $
   $y \cdot (\e^{i-1}, a, \e^{n-i})\big) : \  a \in A_k, \ $
   $(x, y) \in h, \ x_i, y_i \in A_2^{\,*},$

\smallskip

\hspace{3.95in}  and $x, y \in n A_{[1,k[\,} \cdot n A_k^{\,*} \, \}\,$.
\end{minipage}

\medskip

\noindent So $\,{\sf restr}^{A_2 A_k}_i(.)\,$ is not a uniform restriction,
but a restriction defined so as to complete the commutation relation 
$(\star)$. 

\bigskip

\noindent Based on the formulas for ${\sf restr}^{A_2}_i(g)\,$ (Def.\ 
\ref{DEFunifrestrMorph}) and $\iota(.)$, we have:

\bigskip

\noindent $\iota({\sf restr}^{A_2}_i(g))$ 
$ \ = \ $
$\iota \big(\{((\ldots, x_{i-1}, x_i a_0, x_{i+1}, \ldots),\, $
  $(\ldots, y_{i-1}, y_i a_0, y_{i+1}, \ldots)) :\, (x,y) \in g\}$

\smallskip

\hspace{0.9in} $\,\cup\,$
$\{((\ldots, x_{i-1}, x_i a_1, x_{i+1}, \ldots),\, $
   $(\ldots, y_{i-1}, y_i a_1, y_{i+1}, \ldots)) :\, (x,y) \in g\} \big)$

\smallskip

\begin{minipage}{\textwidth}
$= \ {\sf id}_{n A_k \minus n A_{[1,k[\,}}$

\smallskip

\hspace{0.2in} $\cup$
 \ $\{ \big( (u_1,\ldots,u_{i-1}, a_1 x_i^{(r)}a_0, u_{i+1},\ldots,u_n),\,$
$(v_1,\ldots,v_{i-1}, a_1 y_i^{(r)}a_0, v_{i+1},\ldots,v_n) \big)\, : \ $
$1 \le r \le \ell$,

\hspace{0.5in} 
$(\forall s \in [1,n] \minus \{i\})(\exists l \in [2,k[\,)$ 
$[(u_s, v_s)$ $\in$ $\{(a_1x_s^{(r)},\, a_1 y_s^{(r)}),$
   $((a_1x_s^{(r)})'_l,\, (a_1y_s^{(r)})'_l)\}]\}$  

\smallskip

\hspace{0.2in} $\cup$
 \ $\{ \big( (u_1,\ldots,u_{i-1}, a_1 x_i^{(r)}a_1, u_{i+1},\ldots,u_n),\,$
$(v_1,\ldots,v_{i-1}, a_1 y_i^{(r)}a_1, v_{i+1},\ldots,v_n) \big)\, : \,$
$1 \le r \le \ell$,

\hspace{0.5in} 
$(\forall s \in [1,n] \minus \{i\})(\exists l \in [2,k[\,)$
$[(u_s, v_s)$ $\in$ $\{(a_1x_s^{(r)},\, a_1 y_s^{(r)}),$
   $((a_1x_s^{(r)})'_l,\, (a_1y_s^{(r)})'_l)\}]\}$ 

\smallskip

\hspace{0.2in} $\cup$
 \ $\{\big( (u_1,\ldots,u_{i-1},(a_1x_i^{(r)}a_0)'_j,u_{i+1},\ldots,u_n),\,$
$(v_1,\ldots,v_{i-1},(a_1y_i^{(r)}a_0)'_j, v_{i+1},\ldots,v_n) \big)\, : \,$
$1 \le r \le \ell$,

\hspace{0.5in} 
$(\forall s \in [1,n] \minus \{i\})(\exists l \in [2,k[\,)$
$[(u_s, v_s)$ $\in$ $\{(a_1x_s^{(r)},\, a_1 y_s^{(r)}),$
   $((a_1x_s^{(r)})'_l,\, (a_1y_s^{(r)})'_l)\}]\}$

\smallskip

\hspace{0.2in} $\cup$
 \ $\{\big( (u_1,\ldots,u_{i-1},(a_1x_i^{(r)}a_1)'_j,u_{i+1},\ldots,u_n),\,$
$(v_1,\ldots,v_{i-1},(a_1y_i^{(r)}a_1)'_j, v_{i+1},\ldots,v_n) \big)\, : \,$
$1 \le r \le \ell$,

\hspace{0.5in} 
$(\forall s \in [1,n] \minus \{i\})(\exists l \in [2,k[\,)$
$[(u_s, v_s)$ $\in$ $\{(a_1x_s^{(r)},\, a_1 y_s^{(r)}),$
   $((a_1x_s^{(r)})'_l,\, (a_1y_s^{(r)})'_l)\}]\}$.
\end{minipage}

\bigskip

\noindent By \cite[Lemma 2.8]{BiEmbHTarx} (quoted above) this becomes

\medskip

\begin{minipage}{\textwidth}
\noindent $\iota({\sf restr}^{A_2}_i(g))$

\smallskip

$= \ {\sf id}_{n A_k \minus n A_{[1,k[\,}}$

\smallskip

\hspace{0.2in} $\cup$
 \ $\{ \big( (u_1,\ldots,u_{i-1}, a_1 x_i^{(r)}a_0, u_{i+1},\ldots,u_n),\,$
$(v_1,\ldots,v_{i-1}, a_1 y_i^{(r)}a_0, v_{i+1},\ldots,v_n) \big)\, : \,$
$1 \le r \le \ell$,

\hspace{0.5in} 
$(\forall s \in [1,n] \minus \{i\})(\exists l \in [2,k[\,)$
$[(u_s, v_s)$ $\in$ $\{(a_1x_s^{(r)},\, a_1 y_s^{(r)}),$
   $((a_1x_s^{(r)})'_l,\, (a_1y_s^{(r)})'_l)\}]\}$

\smallskip

\hspace{0.2in} $\cup$
 \ $\{ \big( (u_1,\ldots,u_{i-1}, a_1 x_i^{(r)}a_1, u_{i+1},\ldots,u_n),\,$
$(v_1,\ldots,v_{i-1}, a_1 y_i^{(r)}a_1, v_{i+1},\ldots,v_n) \big)\, : \,$
$1 \le r \le \ell$,

\hspace{0.5in} 
$(\forall s \in [1,n] \minus \{i\})(\exists l \in [2,k[\,)$
$[(u_s, v_s)$ $\in$ $\{(a_1x_s^{(r)},\, a_1 y_s^{(r)}),$
   $((a_1x_s^{(r)})'_l,\, (a_1y_s^{(r)})'_l)\}]\}$

\smallskip

\hspace{0.2in} $\cup$
 \ $\{\big( (u_1,\ldots,u_{i-1},(a_1x_i^{(r)})'_j,u_{i+1},\ldots,u_n),\,$
$(v_1,\ldots,v_{i-1},(a_1y_i^{(r)})'_j, v_{i+1},\ldots,v_n) \big)\, : \,$
$1 \le r \le \ell$,

\hspace{0.5in} 
$(\forall s \in [1,n] \minus \{i\})(\exists l \in [2,k[\,)$
$[(u_s, v_s)$ $\in$ $\{(a_1x_s^{(r)},\, a_1 y_s^{(r)}),$
   $((a_1x_s^{(r)})'_l,\, (a_1y_s^{(r)})'_l)\}]\}$

\smallskip

\hspace{0.2in} $\cup$
 \ $\{\big( (u_1,\ldots,u_{i-1}, a_1x_i^{(r)} a_j,u_{i+1},\ldots,u_n),\,$
$(v_1,\ldots,v_{i-1}, a_1y_i^{(r)} a_j, v_{i+1},\ldots,v_n) \big)\, : \,$
$1 \le r \le \ell$,

\hspace{0.5in} 
$(\forall s \in [1,n] \minus \{i\})(\exists l \in [2,k[\,)$
$[(u_s, v_s)$ $\in$ $\{(a_1x_s^{(r)},\, a_1 y_s^{(r)}),$
   $((a_1x_s^{(r)})'_l,\, (a_1y_s^{(r)})'_l)\}]\}\,$.
\end{minipage}

\bigskip

\noindent On the other hand,

\medskip

\begin{minipage}{\textwidth}

\noindent ${\sf restr}_i^{A_2 A_k}(\iota(g))$

\medskip

$= \ $  
${\sf restr}_i^{A_2 A_k} \big($
${\sf id}_{n A_k \minus n A_{[1,k[\,}}$

\smallskip

\hspace{0.2in} $\cup$
 \ $\{ \big( (u_1,\ldots,u_{i-1}, a_1 x_i^{(r)}, u_{i+1},\ldots,u_n),\,$
$(v_1,\ldots,v_{i-1}, a_1 y_i^{(r)}, v_{i+1},\ldots,v_n) \big) : \,$
$1 \le r \le \ell$,

\hspace{0.5in} 
$(\forall s \in [1,n] \minus \{i\})(\exists l \in [2,k[\,)$
$[(u_s, v_s)$ $\in$ $\{(a_1x_s^{(r)},\, a_1 y_s^{(r)}),$
   $((a_1x_s^{(r)})'_l,\, (a_1y_s^{(r)})'_l)\}]\}$

\medskip

\hspace{0.2in} $\cup$
 \ $\{\big( (u_1,\ldots,u_{i-1},(a_1x_i^{(r)})'_j,u_{i+1},\ldots,u_n),\,$
$(v_1,\ldots,v_{i-1},(a_1y_i^{(r)})'_j, v_{i+1},\ldots,v_n) \big) : \,$
$1 \le r \le \ell$,

\hspace{0.5in} 
$(\forall s \in [1,n] \minus \{i\})(\exists l \in [2,k[\,)$
$[(u_s, v_s)$ $\in$ $\{(a_1x_s^{(r)},\, a_1 y_s^{(r)}),$
   $((a_1x_s^{(r)})'_l,\, (a_1y_s^{(r)})'_l)\}]\}$
$\big)$
\end{minipage}

\medskip

\begin{minipage}{\textwidth}
$= \ {\sf id}_{n A_k \minus n A_{[1,k[\,}}$

\smallskip

\hspace{0.2in} $\cup$
 \ $\{ \big( (u_1,\ldots,u_{i-1}, a_1 x_i^{(r)} a, u_{i+1},\ldots,u_n),\,$

\hspace{0.56in} 
$(v_1,\ldots,v_{i-1}, a_1 y_i^{(r)} a, v_{i+1},\ldots,v_n) \big) : \,$
$a \in A_k,$  $1 \le r \le \ell$,

\hspace{0.5in} 
$(\forall s \in [1,n] \minus \{i\})(\exists l \in [2,k[\,)$
$[(u_s, v_s)$ $\in$ $\{(a_1x_s^{(r)},\, a_1 y_s^{(r)}),$
   $((a_1x_s^{(r)})'_l,\, (a_1y_s^{(r)})'_l)\}]\}$

\smallskip

\hspace{0.2in} $\cup$
 \ $\{\big( (u_1,\ldots,u_{i-1},(a_1x_i^{(r)})'_j,u_{i+1},\ldots,u_n),\,$
$(v_1,\ldots,v_{i-1},(a_1y_i^{(r)})'_j, v_{i+1},\ldots,v_n) \big) : \,$
$1 \le r \le \ell$,

\hspace{0.5in} 
$(\forall s \in [1,n] \minus \{i\})(\exists l \in [2,k[\,)$
$[(u_s, v_s)$ $\in$ $\{(a_1x_s^{(r)},\, a_1 y_s^{(r)}),$
   $((a_1x_s^{(r)})'_l,\, (a_1y_s^{(r)})'_l)\}]\}\,$.
\end{minipage}

\bigskip

\noindent Let us show that 
 \ IR $=_{\rm def}$  $\iota({\sf restr}^{A_2}_i(g))$ $=$ 
${\sf restr}_i^{A_2 A_k}(\iota(g))$ $=_{\rm def}$ RI. 
 \ The ${\sf id}_{n A_k \minus n A_{[1,k[\,}}$ parts are identical in IR 
and RI, so we can ignore those from now on.

\smallskip

\noindent [IR $\subseteq$ RI] \ The 1st and 2nd rows of IR are contained in
the 1st row of RI, when $a \in \{a_0, a_1\}$ ($\,\subseteq A_k$).

The 3rd row of IR is equal to the 2nd row of RI.

The 4th row of IR is contained in the 1st row of RI, when 
$a \in A_{[2,k[\,}$ ($\,\subseteq A_k$).

\smallskip

\noindent [IR $\supseteq$ RI] \ The 1st row of RI is contained in the 1st,
2nd, and 4th rows of IR (for respectively $a = a_0$, $a = a_1$, 
$a \in A_{[2,k[\,}$).

The 2nd row of RI is equal to the 3rd row of IR. 

\medskip

\noindent This proves the commutation relation $(\star)$. 

\medskip

The main consequence of the commutativity relation is that if $g'$ is 
obtained from $g$ by coordinatewise restrictions then $\iota(g')$ can be 
obtained from $\iota(g)$ by restrictions, hence 
$\iota(g')$  $\equiv_{\rm end}$ $\iota(g)$.

\bigskip

\noindent (2) To complete the proof that $\iota$ is a homomorphism,
consider $h, g \in$ $n G_{2,1}^{\rm unif}$ with 

\smallskip

 \ ${\rm domC}(g)$ $\,=\,$ $\{x^{(r)} : 1 \le r \le \ell\}$ $\,=\,$ 
{\large \sf X}$_{i=1}^n P_i$, 

\smallskip

 \ ${\rm imC}(g)$ $\,=\,$ ${\rm domC}(h)$ $\,=\,$ 
$\{y^{(r)} : 1 \le r \le \ell\}$ $\,=\,$ {\large \sf X}$_{i=1}^n Q_i$,

\smallskip

 \ ${\rm imC}(h)$ $\,=\,$ $\{z^{(r)} : 1 \le r \le \ell\}$ $\,=\,$
{\large \sf X}$_{i=1}^n R_i$,

\smallskip

\noindent where $P_i$, $Q_i$, and $R_i$ are finite maximal prefix codes in
$A_2^{\,*}$.  We can apply uniform restrictions to $g$ and $h$, without 
changing the elements of $n G_{2,1}^{\rm unif}$ that $g$ and $h$ represent, 
so as to make the image of $g$ equal to the domain of $h$; by Lemma 
\ref{LEMunifprodPres}, the restricted right-ideal morphisms are also 
uniform. 

Let the tables of $h$ and $g$ be 
$\,g = \{(x^{(r)}, y^{(r)}) : 1 \le r \le \ell\}$, and
$\,h = \{(y^{(r)}, z^{(r)}) : 1 \le r \le \ell\}$; so 
$h \circ g(.) = \{(x^{(r)}, z^{(r)}) : 1 \le r \le \ell\}$.
Then: 

\medskip

${\rm domC}(\iota(g))$ $\,=\,$ 
$(n A_k \minus n A_{[1,k[\,})$ \ \ $\cup$
 \ \ {\large \sf X}$_{i=1}^n \big( a_1 P_i$  $\,\cup\,$
 ${\sf spref}(a_1 P_i) \cdot A_{[2,k[\,}\big)\,$,

\medskip

${\rm imC}(\iota(g))$ $\,=\,$ ${\rm domC}(\iota(h))$  $\,=\,$
$(n A_k \minus n A_{[1,k[\,})$ \ \ $\cup$
 \ \ {\large \sf X}$_{i=1}^n \big( a_1 Q_i $
   $\,\cup\,$  ${\sf spref}(a_1 Q_i) \cdot A_{[2,k[\,}\big)\,$,

\medskip

${\rm imC}(\iota(h))$  $\,=\,$
$(n A_k \minus n A_{[1,k[\,})$ \ \ $\cup$
 \ \ {\large \sf X}$_{i=1}^n \big( a_1 R_i $
   $\,\cup\,$  ${\sf spref}(a_1 R_i) \cdot A_{[2,k[\,}\big)\,$,

\medskip

\noindent with the tables of $\iota(g)$ and $\iota(h)$ given by

\medskip

\begin{minipage}{\textwidth}
$\iota(g) \ = \ {\sf id}_{n A_k \minus n A_{[1,k[\,}}$

\smallskip

\hspace{0.5in} $\cup$
 \ $\{  (u,v) \, : \,$
$ \ u \in {\rm domC}(\iota(g)) \cap n A_{[1,k[\,} \cdot nA_k^{\,*}$,
$ \ v \in {\rm imC}(\iota(g)) \cap n A_{[1,k[\,} \cdot nA_k^{\,*}$, \ and

\smallskip

\hspace{0.9in} for all $r \in \{1,\ldots,\ell\}$,
$i \in \{1,\ldots,n\}$, and $j \in [2,k[\,$,

\hspace{1.1in} if $\,u_i = a_1 \, x_i^{(r)} \in a_1 P_i\,$ then
   $\,v_i = a_1 \, y_i^{(r)}\,$;

\hspace{1.1in} if $\,u_i = (a_1 \, x_i^{(r)})'_j$  $\,\in\,$
  ${\sf spref}(a_1 \, P_i) \cdot a_j \ $ then
$\,v_i = (a_1 \, y_i^{(r)})'_j\,\}$;  
\end{minipage}

\medskip

\begin{minipage}{\textwidth}
$\iota(h) \ = \ {\sf id}_{n A_k \minus n A_{[1,k[\,}}$

\smallskip

\hspace{0.5in} $\cup$
 \ $\{  (v,w) \, : \,$
$ \ v \in {\rm domC}(\iota(h)) \cap n A_{[1,k[\,}$,
$ \ w \in {\rm imC}(\iota(h)) \cap n A_{[1,k[\,}$, \ and

\smallskip

\hspace{0.9in} for all $r \in \{1,\ldots,\ell\}$,
$i \in \{1,\ldots,n\}$, and $j \in [2,k[\,$,

\hspace{1.1in} if $\,v_i = a_1 \, y_i^{(r)} \in a_1 Q_i\,$ then
   $\,w_i = a_1 \, z_i^{(r)}\,$;

\hspace{1.1in} if $\,v_i = (a_1 \, y_i^{(r)})'_j$  $\,\in\,$
  ${\sf spref}(a_1 \, Q_i) \cdot a_j \ $ then
$\,w_i = (a_1 \, z_i^{(r)})'_j\,\}$.
\end{minipage}

\medskip

\noindent Similarly, from the table of $h \circ g$ we obtain

\medskip

\begin{minipage}{\textwidth}
$\iota(h \circ g) \ = \ {\sf id}_{n A_k \minus n A_{[1,k[\,}}$

\smallskip

\hspace{0.5in} $\cup$
 \ $\{ (u,w) \, : \,$
$ \ u \in {\rm domC}(\iota(g)) \cap n A_{[1,k[\,}^{\,*}$,
$ \ w \in {\rm imC}(\iota(h)) \cap n A_{[1,k[\,}^{\,*}$, \ and

\smallskip

\hspace{0.9in} for all $r \in \{1,\ldots,\ell\}$,
$i \in \{1,\ldots,n\}$, and $j \in [2,k[\,$,

\hspace{1.1in} if $\,u_i = a_1 \, x_i^{(r)} \in a_1 Q_i\,$ then
   $\,w_i = a_1 \, z_i^{(r)}\,$;

\hspace{1.1in} if $\,u_i = (a_1 \, x_i^{(r)})'_j$  $\,\in\,$
  ${\sf spref}(a_1 \, Q_i) \cdot a_j \ $ then
$\,w_i = (a_1 \, z_i^{(r)})'_j\,$
\}.
\end{minipage}

\medskip

\noindent One observes immediately that $\,\iota(h \circ g)$ $=$
$\iota(h) \circ \iota(g)$.
 \ \ \   \ \ \ $\Box$

\section{The word problem of {\boldmath $n G_{k,1}$} } 

The following is well known.

\begin{lem} \label{LEMcomplexityRel}
 \ Let $M_1$ and $M_2$ be two finitely generated monoids such that
$M_2 \subseteq M_1$.

\smallskip

\noindent {\small \rm (1)} \ If the word problem of $M_1$ is in {\sf coNP}
(or in {\sf P} or in {\sf NP}), then the word problem of $M_2$ is in
{\sf coNP} too (respectively in {\sf P} or in {\sf NP}).

\smallskip

\noindent {\small \rm (2)} \ If the word problem of $M_2$ is {\sf coNP}-hard
or {\sf NP}-hard (with respect to polynomial-time many-one reductions), then 
the word problem of $M_1$ is also {\sf coNP}-hard, respectively 
{\sf NP}-hard.  
\end{lem}
{\sf Proof.} Let $\Gamma_j$ be a finite generating set of $M_j$ ($j = 1,2$).
Since $M_2 \subseteq M_1$, for every $\gamma \in \Gamma_2$ there exists
$w_{\gamma}$ $\in$  $\Gamma_1^{\,*}$ such that $\gamma =_{M_2} w_{\gamma}$.
For $\,x = \gamma_1 \ \ldots \ \gamma_m$ $\in \Gamma_2^{\,*}$, let
$w_x \in \Gamma_1^{\,*}$ be the concatenation
$\,w_{\gamma_1} \ \ldots \ w_{\gamma_m}$.

\smallskip

\noindent {\small \rm (1)} Let $\cal A$ be a (co-non)deterministic algorithm
that on input $(x, y) \in \Gamma_1^{\,*} \x \Gamma_1^{\,*}$ decides whether
$x =_{M_1} y$.  Let $(u,v) \in \Gamma_2^{\,*} \x \Gamma_2^{\,*}$;
then $\cal A$, applied to $(w_u, w_v) \in \Gamma_1^{\,*} \x \Gamma_1^{\,*}$
decides whether $u =_{M_2} v$ (since the latter is equivalent to
$w_u =_{M_1} w_v$.

\smallskip

\noindent {\small \rm (2)} Let $L \subseteq A^*$ be any language in
{\sf coNP}, and let $\,\rho: x \in A^* \to (\rho_1(x),\, \rho_2(x))$ $\in$
$\Gamma_2^{\,*}$ be a polynomial-time many-one reduction of $L$ to the word
problem of $M_2$ (i.e., $x \in L$ iff $\rho_1(x) =_{M_2} \rho_2(x)$).
Then $\,x \in A^*$ $\longmapsto$ $(\rho_1(x),\, \rho_2(x))$
$\longmapsto$  $(w_{\rho_1(x)},\, w_{\rho_2(x)})\,$
is a polynomial-time many-one reduction of $L$ to the word problem of $M_1$.
 \ \ \ $\Box$

\begin{lem} \label{LEMwpSigmaTauF}
 \ The word problem of the finitely generated subgroup 
$\,\langle \sigma,$ $\tau_{2,1} \x {\mathbb 1},$ 
{\small \sf F} $\x$ ${\mathbb 1}\rangle\,$ of $\,2 G_{2,1}^{\rm unif}$
is {\sf coNP}-hard (with respect to polynomial-time many-one reduction).
\end{lem}
{\sc Proof.} In \cite[Thm.\ 4.22]{BinG} we proved that the word problem of 
the infinitely generated subgroup 
$\,\langle \tau \,\cup\, \{{\small \sf F}\} \rangle\,$ of $G_{2,1}$ is 
{\sf coNP}-hard.
This subgroup represents bijective circuits, whose equivalence problem was
proved {\sf coNP}-complete by Stephen Jordan \cite{Jordan}.

In \cite[Lemma 4.23]{BinG} we embedded 
$\,\langle \tau \,\cup\, \{{\small \sf F}\} \rangle_{G_{2,1}}\,$ first into 
$\,\langle \tau \x \{{\mathbb 1}\}$ $ \,\cup \, $
$\{{\small \sf F} \x {\mathbb 1}\} \rangle_{2 G_{2,1}}\,$, and then into
$\,\langle \sigma,$ $\tau_{2,1} \x {\mathbb 1},$ 
{\small \sf F} $\x$ ${\mathbb 1} \rangle_{2 G_{2,1}}\,$ by 

\medskip

$\tau_{j,j+1} \ \mapsto \ \sigma^{j-1}$ $\circ$ 
$(\tau_{2,1} \x {\mathbb 1})$ $\circ$  $\sigma^{-j+1}$
 \ \ \ (for all $\tau_{j,j+1} \in \tau$).

\medskip
 
\noindent This shows that the word problem of $2 G_{2,1}$ over any finite
generating set is {\sf coNP}-hard; it also shows that the word problem of
the finitely generated subgroup
$\,\langle \sigma,$ $\tau_{2,1} \x {\mathbb 1},$ 
{\small \sf F} $\x$ ${\mathbb 1}\rangle\,$ of $2 G_{2,1}$ is 
{\sf coNP}-hard.

And we saw in Lemma \ref{LEMexANdcouterex} that $\,\sigma,$ 
$\tau_{2,1} \x {\mathbb 1},$ {\small \sf F} $\x$ ${\mathbb 1}$ $\in$ 
$2 G_{2,1}^{\rm unif}$. 
 \ \ \ $\Box$

\begin{thm} \label{THMwp}
 \ The word problem of $n G_{k,1}$ over a finite generating set is 
{\sf coNP}-complete (with respect to polynomial-time many-one reduction), 
 for all $n \ge 2$ and $k \ge 2$.
\end{thm}
{\sc Proof.} By Theorem \ref{THMnGk1FinGen}, $n G_{k,1}$ is finitely 
generated. 
In \cite[Lemma 3.5]{BinG} it was proved that the word problem 
of $n G_{2,1}$ is in {\sf coNP}.  The proof also works for $n G_{k,1}$, and
actually does not depend on $k$, but only on the fact that $n G_{k,1}$ has a
finite generating set.
Hence the word problem of $n G_{k,1}$ is in {\sf coNP}.

\smallskip

By Lemma \ref{LEMwpSigmaTauF}, the word problem of the subgroup of 
$2 G_{2,1}^{\rm unif}$ generated by the finite set $\,\{\sigma,$ 
$\tau_{2,1} \x {\mathbb 1},$ {\small \sf F} $\x$  ${\mathbb 1}\}\,$ is 
{\sf coNP}-hard. Since $2 G_{2,1}^{\rm unif}$ embeds into $2 G_{k,1}$,
the word problem of $2 G_{k,1}$ is {\sf coNP}-hard (by Lemma 
\ref{LEMcomplexityRel}). 
And $2 G_{k,1}$ obviously embeds into $n G_{k,1}$.  Hence the problem of 
$n G_{k,1}$ over a finite generating set is {\sf coNP}-hard.
 \ \ \ $\Box$

\bigskip

\bigskip




\medskip

\noindent {\scriptsize  birget@camden.rutgers.edu }

\end{document}